\documentclass[a4paper, 11pt]{article} 
\usepackage{amsfonts,amsmath,amssymb,amscd}
\usepackage[all]{xy} 
\usepackage{latexsym}
\usepackage[T1]{fontenc}
 
\newtheorem{thm}{Theorem}[section]
\newtheorem{definition}[thm]{Definition}
\newtheorem{prp}[thm]{Proposition}
\newtheorem{lemma}[thm]{Lemma}
\newtheorem{cor}[thm]{Corollary}
\newtheorem{moncor}[thm]{Corollary}
\newtheorem{remark}[thm]{Remark}

\newenvironment{prf}[1][]{\begin{trivlist}\item[]{\bf Proof#1. }}
{\hfill $\blacksquare$ \end{trivlist}}

\newcommand{\nap}{\nabla^\perp}
\newcommand{\coker}{\mathrm{coker }}
\newcommand{\ind}{\mathrm{index}\,}

\newcommand{\bc}{\mathcal{C}}
\newcommand{\bd}{\mathcal{D}}
\newcommand{\grad}{\vec{\nabla}}  
\newcommand{\bfff}{\mathcal{F}}
\newcommand{\beee}{\mathcal{E}}
\newcommand{\bp}{\mathcal{P}}
\newcommand{\br}{\mathcal{R}}
\newcommand{\bl}{\mathcal{L}}
\newcommand{\bs}{\mathcal{S}}

\newcommand{\bm}{\mathcal{M}}

\newcommand{\ba}{\mathcal {A}}

\newcommand{\bz}{\mathcal {Z}}

\newcommand{\na}{\nabla}
\newcommand{\be}{\begin{equation}}
\newcommand{\ee}{\end{equation}}

\newcommand{\si}{\sigma}

\newcommand{\kah}{\text{K\"{a}hler }}

\newcommand{\beq}{\begin{eqnarray*}}
\newcommand{\eeq}{\end{eqnarray*}}
\newcommand{\bpr}{\begin{prf}}
\newcommand{\epr}{\end{prf}}

\newcommand{\Z}{\mathbb{Z}}
\newcommand{\T}{\mathbb{T}}

\newcommand{\R}{\mathbb{R}}

\newcommand{\C}{\mathbb{C}}

\newcommand{\dsurdt}{\frac{\partial}{\partial t}}
\newcommand{\half}{\frac{1}{2}}
\newcommand{\quart}{\frac{1}{4}}
\newcommand{\troisquart}{\frac{3}{4}}

\begin{document}
\title{Smooth moduli spaces of associative submanifolds}
\author{Damien Gayet}
\maketitle

\centerline{\textbf{Abstract}}
Let $M^7$ be a smooth manifold equipped with a $G_2$-structure $\phi$,
and $Y^3$ be a closed compact  $\phi$-associative submanifold.  In \cite{McL}, R. McLean proved that the moduli space  $\bm_{Y,\phi}$ of the $\phi$-associative deformations of $Y$ has vanishing virtual dimension. 
In this paper, we perturb  $\phi$ into a  $G_2$-structure $\psi$ in order to ensure  the smoothness of $\bm_{Y,\psi}$ near $Y$. If $Y$ is allowed to have a boundary moving in a fixed coassociative submanifold $X$, 
it was proved in \cite{GaWi} that the moduli space $\bm_{Y,X}$ of the associative deformations of $Y$ with boundary in $X$ has finite virtual dimension. We show here that a generic perturbation of the boundary condition $X$ into $X'$ gives the smoothness of  $\bm_{Y,X'}$. In another direction, we use  Bochner's technique to prove a vanishing theorem that forces $\bm_Y$ or $\bm_{Y,X}$ to be  smooth near $Y$. For every case, some explicit families of examples will be given. 
\bigskip


\smallskip

\section{Introduction}
In the Euclidean space $(\R^7,g_0)$ with its canonical coordinates $(x_i)_{i=1, \cdots, 7}$, consider the 3-form 
$$\phi_0 = dx_{123}+ dx_{145}+ dx_{167}+ dx_{246}- dx_{257} - dx_{347}-dx_{356}$$
and $G_2$ the subgroup of $SO(7)$ defined by $G_2 =  \{g\in SO(7), g^*\phi_0 = \phi_0\}$. If $M$ is an oriented spin 7-dimensional Riemannian manifold, its structure group can be reduced to  $G_2\subset SO(7)$. Given a set of trivialization charts for $TM$ compatible with $G_2$, $M$ inherits a nondegenerate 3-form $\phi$ and a metric $g$, which
are the pullbacks of $\phi_0$ and $g_0$ by these charts. 
We call the pair $(\phi,g)$ a $G_2$-$structure$. 
 Moreover, $TM$ inherits a vector product $\times$ defined by $$\forall u,v,w \in TM, \langle u\times v, w\rangle  = g(u\times v,w) = \phi(u,v,w).$$
Note that in $\R^7$, the subspace $\R^3\times \{0\}$ is stable under this vector product, which induces the classical vector product on $\R^3$.  When  $\phi$ is closed and coclosed for $g$, the structure is said to be \textit{torsion-free}. In this situation, 
 the holonomy of $g$ is a subgroup of $G_2$, see \cite{Joyce}.   \\
 
\noindent
A 3-dimensional submanifold $Y$ in $(M,\phi,g)$ is called $\phi$-\textit{associative}, or simply  \textit{associative} when
there is no ambiguity, if its tangent bundle is stable under the
vector product associated to $\phi$. In other terms, 
$\phi$ restricted to $Y$ is a volume form for $Y$.
Likewise, a 4-dimensional submanifold $X$ 
is called \textit{coassociative} if  the fibers of its
normal bundle are associative, or equivalently, 
$\phi_{|TX}$ vanishes. 


\subsection{Genericity}
\textbf{Closed associative submanifolds.} 
\begin{definition} 
Consider a smooth spin 7-manifold $M$ and $Y$ a smooth compact closed 3-submanifold. For every $G_2$-structure $\phi$, define
$ \bm_{Y,\phi}$ to be the set of smooth $\phi$-associative
submanifolds isotopic to $Y$. 
\end{definition}
It is known from \cite{McL} that the problem of associative deformations of a compact closed associative submanifold $Y$  is related to an elliptic partial differential equation, namely a twisted Dirac operator, see Theorem \ref{Dirac}. Hence for a fixed $G_2$-structure $\phi$, the moduli space
$\bm_{Y,\phi}$ has finite and vanishing virtual  dimension. In general, the situation is obstructed. For instance,
consider  the  torus $\T^3\times \{t\}$ in the flat torus $(\T^7,\phi_0,g_0) = \T^3\times \T^4$. 
This is an associative submanifold, and its moduli space $\bm_{\T^3\times \{t\}}$ of associative deformations contains at least the 4-dimensional $\T^4$.  See also Proposition \ref{Calabi} for a more general situation in a product of a Calabi-Yau manifold with 
$S^1$. \\

\noindent
A natural question is to find conditions which force the moduli space $\bm_{Y,\phi}$ to be smooth at least near a $\phi$-associative $Y$, or in other terms, which force the cokernel of the operator to vanish. 
One  way to solve this  is to perturb the $G_2$-structure
and get generic smoothness. It turns out that in general we cannot do this in the realm of  torsion-free  structures, see Remark \ref{torsion-free}. On the other hand, 
$G_2$-structures with \textit{closed} 3-form $\phi$ seem to be rich enough to work with, at least for the point of view of calibrated geometries, see \cite{HaLa}. Indeed, any closed $G_2$-structure $\phi$ defines a calibration, and when this form is closed, the calibrated submanifolds, here the associative ones, do minimize the volume in their homology class. As suggested to the author by D. Joyce, we will prove the following 
\begin{thm} \label{generic} \textit{Let $M$ be a manifold equipped with a closed $G_2$-structure $\phi$, and $Y$ be a smooth  compact closed $\phi$-associative submanifold. Then there is a neighbourhood $V$ of $Y$, such that for every generic closed $G_2$-structure $\psi$ close enough to $\phi$, the subset of elements of $\bm_{Y,\psi}$ lying in $V$ is a finite set, possibly empty.}
\end{thm}
A former result in this direction was proved by S. Abkulut and S. Salur \cite{AkSa2}, where the authors allow a certain freedom for the definition of associativity. \\

\noindent
\textbf{Associative submanifolds with boundary.}
In \cite{GaWi}, the authors showed that 
the problem of associative deformations of an associative submanifold $Y$ with boundary in 
a fixed coassociative submanifold $X$ is an elliptic problem of finite index. Moreover, they proved that this virtual dimension equals the index of a natural Cauchy-Riemann operator related to the complex geometry of the boundary, see Theorem \ref{boundary} below. 
As in the case of a closed associative,
the situation can be obstructed. For instance, consider in $(\T^7, \phi_0, g_0)$
the $T^2$- family of associative submanifolds $$Y_\lambda = 
\{(x_1,x_2,x_3, \lambda,\mu,0,0), 0\leq x_1\leq 1/2, x_2, x_3 \in S^1\}, (\lambda, \mu) \in T^2.$$  
The two components of the boundary of $Y_\lambda$ lie
in the union $X$ of the two  coassociatives tori $$X_i = \{(i/2, x_2, x_3, x_4, x_5,0,0), x_2, x_3, x_4, x_5 \in S^1\}, i=0,1.$$ However
the index of this problem vanishes, see \cite{GaWi} or Theorem \ref{boundary}.  For  more general obstructed situations, see Theorem \ref{SL-boundary}.\\

\noindent
As in the case of a closed associative, we can perturb the closed $G_2$-structure  $\phi$ of the manifold $M$ into $\psi$  to ensure the smoothness of the moduli space. 
Note that in this case, $X$ has no reason to remain coassociative for the new structure. But it remains $\psi$-\textit{free}, i.e the tangent
space of $X$  does not contain any $\psi$-associative 3-plane, see \cite{HaLa2} or \cite[Section 5]{GaWi}. Indeed, $\phi$-coassociativity implies $\phi$-freedom, and for a submanifold to be  $\phi$-free is an open condition in the variable $\phi$. For any $G_2$-structure $\phi$, the problem of deformations of an associative submanifold
with boundary in a fixed $\phi$-free submanifold is still elliptic \cite{GaWi} and, in our present case, its index is the same as the index for the unperturbed situation.
\begin{definition} 
Consider a manifold equipped with a $G_2$-structure $(\phi,g)$ and $Y$ a smooth compact  associative submanifold with boundary in a $\phi$-free submanifold $X$.  We denote by 
$ \bm_{Y,X}$ the set of smooth associative
submanifolds with boundary in $X$ and isotopic to $Y$. 
\end{definition}
Instead of changing the $G_2$-structure, we can move 
the boundary condition, namely $X$. Still, if we demand that $X$ remains coassociative, in general we can not get smoothness. Indeed, it is known \cite{McL} that the moduli space of coassociative perturbations of $X$ is smooth and has the  dimension $b_2^+(X)$ of the space of harmonic self-dual 2-forms on $X$. In the former example of the flat torus, every coassociative deformation of $X$ is a translation of the initial situation, hence the problem  remains obstructed. Now, since any perturbation of a  $\phi$-free submanifold remains $\phi$-free,  we can fix $\phi$ and  perturb $X$.  
\begin{thm}\label{coasso} Let $Y$ be a smooth associative submanifold
with boundary in a  smooth coassociative submanifold 
$X$. If the virtual dimension of $\bm_{Y,X}$ is non-negative, then for any sufficiently small generic smooth
deformation $X'$ of $X$, either $\bm_{Y,X'}$ is locally empty, that means there is no 
associative manifold with boundary in $X'$ close enough to $Y$, or 
there exists a small associative deformation $Y'$ of $Y$ such that the moduli space $\bm_{Y',X'}$  
is smooth near $Y'$ and of  dimension equal to 
the index computed for the unperturbed situation. 
\end{thm}
 \subsection{Metric conditions}
Concrete examples
 are often non generic, so we would like too to get a condition that is not a perturbative one. For holomorphic curves in dimension 4, there are topological conditions on the degree of the normal bundle which imply the smoothness
of the moduli space of complex deformations, see  
 \cite{HoLiSi}. The main reason  is that holomorphic
 curves intersect positively. In our case, there is no such phenomenon. \\

In \cite[page 30]{McL}, R. McLean gives an example 
of an isolated associative submanifold. For this, he recalls that  R. Bryant and S. Salamon constructed in  \cite{BrSa} a metric of holonomy $G_2$ on the spin bundle $S^3\times \R^4$ of the round 3-sphere. In this case, the base $Y= S^3\times \{0\}$ is 
associative, the normal bundle of $Y$ is  the spin bundle of $S^3$, and the operator related to the associative deformations of $Y$ is  the Dirac operator on $S^3$. By the famous theorem of Lichnerowicz \cite{Li}, there are no non trivial harmonic spinors on $S^3$ for metric reasons (to be precise, because the Riemannian scalar curvature is positive), so the sphere is isolated as an associative submanifold. \\

\noindent 
\textbf{Minimal submanifolds.}
Recall that in a manifold with a closed $G_2$-structure,
associative submanifolds are minimal. In \cite{Si}, J. Simons gives a metric condition for a minimal submanifold
to be stable, i.e.  isolated. For this, he introduces the following operator, a sort
of partial Ricci operator:
\begin{definition} Let $(M,g)$ be a Riemannian manifold and $Y$ a $p$-dimensional submanifold in $M$ and $\nu$ be its normal bundle. Choose $\{e_1,\cdots e_p\}$ 
a local orthonormal frame field of $TY$, and define the 0-order operator 
$\br: \Gamma(Y,\nu) \to \Gamma(Y,\nu)$ with $
\br s = \pi_\nu \sum_{i=1}^p  R (e_i, s) e_i,$
where $R$ is the curvature tensor of $g$ on $M$   and $\pi_\nu$
the orthogonal projection to $\nu$.
\end{definition}
 It turns out that the definition is independent of the chosen oriented orthonormal frame, and that
$\br$ is symmetric. Simons defines another operator $\mathcal A$ related to the second fondamental form of $Y$:
\begin{definition} Let $SY$ be  the bundle over $Y$ whose fiber at a point $y$ is the space of symmetric endomorphisms
of $T_yY$, and $A \in Hom (\nu, SY)$ the second fundamental form defined by 
$  A(s) (u) = -\nabla^\top_u s,$ where $u\in TY$, $s \in \nu$, and $\nabla^\top$ is the projection to $TY$
of the  ambient Levi-Civita connection $\nabla$, with $\nabla =
\nabla^\top + \nabla^\perp$.
Denote by $\mathcal A$ the operator $\mathcal A: \Gamma(Y,\nu) \longrightarrow \Gamma(Y,\nu)$, $
\ba s = A^t\circ A (s),
$
where $A^t$ is the transpose of $A$.
\end{definition}
It is classical that $\mathcal A$  is a symmetric positive 0-th order operator. Moreover, it vanishes if $Y$ is totally geodesic. 
Using both operators and  Bochner's technique, Simons gives a sufficient condition for a minimal submanifold 
to be stable:
\begin{thm}[\cite{Si}]\label{Simons}
Let $Y$ be a minimal submanifold in $M$, and assume that $\br - \mathcal A$ 
is positive. Then $Y$ cannot be deformed as a minimal submanifold. 
\end{thm}
In particular, if $Y$ is a compact closed associative submanifold 
satisfying the conditions of Theorem \ref{Simons} in a manifold $M$ with a closed $G_2$-structure, then 
it cannot be perturbed as an associative submanifold.  
Now, if $Y$ is an associative submanifold with a boundary, we introduce another operator:
\begin{definition}\label{delta} In a manifold equipped with a $G_2$-structure, let $Y$ be a smooth compact
associative submanifold with boundary and $\nu$ be its
normal bundle. Let $L$ be a two dimensional real  subbundle of $\nu_{|\partial Y}$ invariant under the action of $n\times$, where
$n$ is the inward unit normal vector field along $\partial Y$. Choose $\{v,w= n\times v\}$ a local orthonormal frame for $T\partial Y$. We denote by $\bd_L$ the operator 
$\bd_L: \Gamma(\partial Y,L) \to  \Gamma(\partial Y,L)$, 
$$ \bd_L s =   \pi_L(v\times \nabla^\perp_w s- w\times \nabla^\perp_v s),$$
where $\pi_L: \nu_{|\partial Y} \to L$ is the orthogonal projection to $L$ and $\nap$ the normal connection on $\nu$ induced
by the Levi-Civita connection $\nabla$ on $M$.
\end{definition}
\begin{remark}
Note that such subbundles always exist. Firstly, it is easy to check that $\nu_{|\partial Y}$ is stable under the action of $n\times$. Secondly, $\nu_{|\partial Y}$  has real dimension $4>2$, so that it has a non vanishing section $e$. 
Then, $L$ generated by $e$ and $n\times e$ satisfies the conditions of the Definition \ref{delta}.
\end{remark}
We will prove in Proposition \ref{bd} that $\bd_L$ is independent of the chosen oriented frame, is of order 0  and is symmetric. Assume further that the boundary of $Y$ lies  in a coassociative submanifold $X$. It turns out that $Y$ intersects $X$ orthogonally, see Theorem \ref{boundary} below. Denote by $\mu_X$ the 2-dimensional orthogonal complement of $n$ in the normal bundle of $X$ over $\partial Y$, where $n$ is the inward normal unit vector field in $Y$ along $\partial Y$.  Then we can state the following vanishing 
\begin{thm}\label{open}
Let $M$ be a manifold equipped with a torsion-free $G_2$-structure and  $Y$ be an associative submanifold  with boundary in a coassociative submanifold $X$. 
If  $\bd_{\mu_X}$ and 
$\br-\mathcal A $ are positive, the moduli space $\bm_{Y,X}$ is  smooth near $Y$ and of dimension 
given by the index in Theorem \ref{boundary}.
\end{thm}
Thanks to Theorem \ref{open}, we can find an explicit example, in the Bryant-Salamon manifold with $G_2$-holonomy, of a locally
smooth one dimensional moduli space of associative deformations
with boundary in a coassociative submanifold, see Corollary \ref{BrySal}.
In Section \ref{exemples}, we explain other explicit examples, in particular for an ambient manifold which is the product of a Calabi-Yau manifold with $S^1$ or $\R$, see Theorem \ref{SL-boundary}. \\

\noindent
 \textit{Acknowledgements. } The author benefits the support of the French Agence nationale de la recherche. Part of this work was done during a visit at the Poncelet Laboratory in Moscow. I am grateful to this institution for its hospitality. I would like to thank Vincent Borrelli (resp. Jean-Yves Welschinger) who convinced me that there is a life after curvature tensors (resp. Sobolev spaces), Gilles Carron and Alexei Kovalev for their interest in this work and Dominic Joyce for a stimulating discussion.
\section {Closed associative submanifolds}
\subsection{The operator $D$ and the deformation problem}
 We begin with the version of McLean's theorem proposed by Akbulut and Salur, and a proof of it.
\begin{thm}[\cite{McL},\cite{AkSa}]\label{Dirac}\label{theo MacLean} Let  $M$ be a manifold equipped with a $G_2$-structure $(\phi,g)$, and $Y$ be a closed compact associative submanifold with  normal bundle $\nu$. Then the Zariski tangent space at $Y$ of $\bm_Y$ 
can be identified with the kernel of the operator $D:  \Gamma(Y,\nu) \to  \Gamma(Y,\nu),$ where 
\begin{eqnarray}\label{Operator}
Ds =  \sum_{i=1}^3 e_i \times \nabla^\perp_{e_i} s + 
\sum_{k=1}^{4}(\nabla_s {*}\phi) (\eta_k, \omega)\otimes \eta_k.
\end{eqnarray} 
Here $(e_i)_{i = 1, 2, 3} $ is
any local  orthonormal frame of the tangent space of $Y$ with $e_3 = e_1\times e_2$,  
$\omega = e_1\wedge e_2\wedge e_3$, $(\eta_k)_{k = 1, 2, 3,4} $ is
any local orthonormal frame of $\nu$ and
 $\nabla^\perp$ is the connection on $\nu$ induced by the Levi-Civita connection $\nabla$ of $(M,g)$.
\end{thm}
 Note that second part 
is a 0-th order operator that vanishes for a torsion-free 
$G_2$-structure, as proved in \cite{AkSa}.
\bpr Firstly, recall the existence on $(M,\phi,g)$ of an 
important object $\chi$, the 3-form with values in $TM$ and defined, for $u$, $v$, $w\in TM$
 by  
 $\chi(u,v,w) = - u \times (v\times w) - \langle u,v\rangle  w + \langle u,w\rangle   v.$
  It is easy to check  \cite{AkSa} that  $\chi(u,v,w) $
 is orthogonal to the 3-plane $u\wedge v\wedge w$. Moreover we will use the  following  useful formula \cite{HaLa}: 
$$ \forall u,v,w,\eta\in TM, \langle \chi(u,v,w), \eta\rangle  = {*}\phi (u,v,w,\eta),$$
where  ${*}$ is the Hodge star associated to the metric $g$. So 
\begin{eqnarray}\label{star}
\chi  = \sum_k \big(\eta_k \lrcorner {*}\phi \big)\otimes \eta_k,
\end{eqnarray}
where $(\eta_k)_{k = 1, 2, \cdots 7} $ is
an local orthonormal frame of the tangent space of $M$.
Further,
 if $Y$ is a 3-dimensional submanifold in $(M, \phi)$, 
 then $\chi_{|TY} = 0$ if and only if $Y$ is associative. 
As in \cite{McL}, we use this characterization 
to study the moduli space of associative deformations of
an associative $Y$.  
Let $Y$ be any smooth closed associative submanifold in $M$. 
We parametrize its deformations by
the sections of its normal bundle $\nu$. Fix $\omega$  a non vanishing global section of $ \Lambda^3 TY$ writing locally $\omega = e_1\wedge e_2\wedge e_3$, with $(e_i)_{i = 1,2,3}$  a local orthonormal frame of $TY$ satisfying $e_3 = e_1\times e_2$. 
For every smooth section $\sigma \in \Gamma(Y,\nu)$, 
define 
\begin{equation}\label{F}
F(\sigma) = \exp_\sigma^*\chi (\omega) \in  \Gamma (Y,  \nu_\sigma),
\end{equation}
where  $\nu_\sigma$ is the normal bundle of 
$\exp_\sigma (Y)$. 
Then $\exp_\sigma (Y)$ is associative if and only if  $F(\sigma)$ vanishes. In order to compute the Zariski tangent space of $\bm_Y$ at the vanishing section, consider a path of normal sections  $(\sigma_t)_{t\in [0,1]}\in \Gamma(Y,\nu)$
and  
$$s = \frac{d \sigma_t}{d t}_{|t=0} \in \Gamma (Y,  \nu).$$
To differentiate $F$ at $\sigma=0$ in the direction of $s$, we use the Levi-Civita connection of $(M,g)$. We have
\begin{eqnarray}\label{derivee} 
 \nabla_{\dsurdt} F(\sigma_t)_{|t=0} = \sum_k \bl_s (\eta_k\lrcorner {*}\phi) (\omega) \otimes \eta_k + (\eta_k\lrcorner {*}\phi) (\omega) \otimes \na_s \eta_k,
\end{eqnarray} 
where $ \bl_s $ is the Lie derivative in the direction $s$.
Since $Y$ is associative,  $\omega\lrcorner {*}\phi = 0 $ and 
the second term   vanishes. Note that this implies that 
the result does not depend  on the chosen connection. 
Thanks to classical Riemannian formulas, we compute the summand of the first term. For every $k$, 
$$ \bl_s (\eta_k\lrcorner {*}\phi) = \eta_k \lrcorner \ \bl_s ({*}\phi) + [\eta_k,s] \lrcorner {*}\phi, $$
and 
since $([\eta_k,s]\wedge \omega) \lrcorner {*}\phi = 0, $
we get 
\begin{eqnarray}\label{derivee} 
 \nabla_{\dsurdt} F(\sigma_t)_{|t=0} = \sum_k \bl_s ({*}\phi) (\eta_k, \omega)\otimes \eta_k.
\end{eqnarray} 
The Lie derivatives can be expressed in terms of the Levi-Civita connection, see for instance Formula 3.3.26 of \cite{Jost},
so that
\begin{eqnarray*}
 \bl_s ({*}\phi) (\eta_k, \omega) & = & (\nabla_s {*}\phi) (\eta_k, \omega) + {*}\phi (\nabla_{\eta_k}s, \omega) + \\
&& {*}\phi (\eta_k, \nabla_{e_1}s , e_2,e_3)+ {*}\phi (\eta_k, e_1,\nabla_{e_2}s ,e_3) +
{*}\phi (\eta_k, e_1,e_2, \nabla_{e_3}s ).
\end{eqnarray*}
The second term of the right hand side vanishes because $\omega\lrcorner {*}\phi =0$
and the third one equals
$ {*}\phi (\eta_k, \nap_{e_1}s , e_2,e_3) = - \langle \nap_{e_1}s \times ( e_2 \times e_3),\eta_k\rangle  .$
Using  the relation $e_2\times e_3 = e_1$ and adding up the two last similar terms,  
we obtain
$ \nabla_s F =  \sum_i e_i \times \nap_i s + \sum_k (\nabla_s {*}\phi) (\eta_k, \omega)\otimes \eta_k.
$ Since $F(0)$ has values in $\nu$, in fact we can assume that the $\eta_k$'s form a local orthonormal frame of $\nu$. 
\epr
\begin{prp}\label{smoothness}
Let $Y$ be a smooth closed associative submanifold in a manifold $M$ equipped with a $G_2$-structure. If the (co)kernel
of  the operator $D$ given by (\ref{Operator})  vanishes, then 
 $ \bm_Y $ is  smooth near $Y$ and of vanishing dimension. In 
 particular, $Y$ is isolated among associative submanifolds isotopic to $Y$.
 \end{prp}
 \bpr 
Fix $Y$ a smooth closed associative submanifold. 
For $kp> 3$, it makes sense to consider the Banach space 
$\beee = W^{k,p}(Y,\nu) $
of sections with weak  derivatives in
$L^p$, up the $k$-th one.  Moreover for $(k-r)/3 >   1/p$, the inclusion 
$W^{k,p}(Y,\nu) \subset C^{r} (Y,\nu) $ holds 
and so $\sigma \in \beee$ is $C^1$ if $k> 1+3/p$. 
In particular, one can define $\nu_\sigma$ the normal bundle
to $\exp_\sigma (Y)$, and  $\bfff$ the Banach bundle over $\beee$ with fiber $ \mathcal{F}_{\sigma} = W^{k-1,p}(Y,\nu_\sigma).$
It is clear that the operator $F$ defined by ($\ref{F}$)  extends to a section $F_{k,p}$ of $\bfff$ over $\beee$.
The proof of Theorem \ref{Dirac} shows that $F_{k,p}$ is smooth and  the derivative of $F$ in the direction of a vector field $s \in T_0\beee= W^{k,p}(Y,\nu)$ is computed 
by (\ref{Operator}).
Now, the operator $D: \Gamma(Y,\nu) \to \Gamma(Y,\nu)$ has symbol $$\sigma (\xi): s \mapsto \sum_i \xi_i s \times e_i = s\times \xi,$$
which is always invertible on $\nu$ as long as  $\xi \in TY\setminus \{0\}.$
This proves that $D$ is elliptic. 
Note that $\sigma (\xi)^2s = - |\xi|^2 s$, which is the symbol
of the Laplacian. 
Hence $F$ is a Fredholm operator, and $\ker D$ and  $ \coker D$ have finite dimension. 
By the implicit function theorem for Banach bundles, 
if $\coker D =\{0\}$, then $F^{-1} (0)$ is  a  smooth Banach submanifold
of $\beee$ near the null section  and of finite dimension equal to  $\dim \ker D = \ind D$, which vanishes since  $Y$ is odd-dimensional.
Lastly, still thanks to the ellipticity of $D$, all elements of $\bm_Y$ are smooth.
\epr
\subsection{Varying the $G_2$-structure}\label{varying G2}
  \textbf{Theorem \ref{generic} } 
\textit{Let $M$ be a manifold equipped with a closed $G_2$-structure $\phi$, and $Y$ be a smooth  compact closed $\phi$-associative submanifold. Then there is a neighbourhood $V$ of $Y$, such that for every generic closed $G_2$-structure $\psi$ close enough to $\phi$, the subset of elements of $\bm_{Y,\psi}$ lying in $V$ is a finite set, possibly empty.}
\bpr  Consider $Y$ a smooth closed associative submanifold in a manifold $M$ equipped with a closed $G_2$-structure $(\phi,g)$. We modify the former map $F$ defined in (\ref{F})
in the following way. For every normal section 
 $\sigma \in \Gamma(Y,\nu)$ and every
$G_2$-structure $\phi'$,  consider
\begin{equation}\label{F psi}
F(\sigma,\phi') = \exp_\sigma^*\chi_{\phi'} (\omega) \in  \Gamma (Y,  \nu_\sigma).
\end{equation}
Here the exponential map corresponds to the fixed metric $g$, whereas $\nu_\sigma$,  the normal vector bundle over $\exp_\sigma (Y)$, depends now on the metric associated to $\phi'$, as does $\chi_{\phi'}$.
We will differentiate $F(0,.)$ in the direction of  
$\bz^3(M)$, the subspace of smooth closed
3-forms on $M$.  Recall that the set of 3-forms defining a $G_2$-structure  
is open in $\Omega^3(M)$, hence for every $\psi \in \bz^3(M)$ with small enough norm, $\phi + \psi $ still defines a closed
$G_2$-structure. Let $(\phi_t)_{t\in [0,1]}$ be a smooth path
of closed $G_2$-structures, with $\phi_0 = \phi$. In formula (\ref{star}), the local orthonormal trivializations $\eta_k$ of the tangent bundle $TM$ are orthonormal for the metric $g_t$ associated to $\phi_t$, consequently we have to choose them as functions of $t$. On the other hand, we can keep $\omega$ constant. Hence $F(0,\phi_t) =  \sum_k ( \eta_k(t)\wedge \omega) \lrcorner {*}_t\phi_t \otimes \eta_k(t)$, where ${*}_t$ denotes the Hodge star for 
 $g_t$. 
Since $\omega\lrcorner {*}\phi=0$, at $t=0$ the
two terms in the derivative containing  $\nabla_{\dsurdt} \eta_k$ vanish, and we have 
\begin{equation*}
 \nabla_{\dsurdt} F(0, \phi_t)_{|t=0} = \sum_k (\eta_k\wedge \omega)\lrcorner 
\dsurdt\Theta(\phi(t))_{|t=0}\otimes \eta_k.
\end{equation*}
The nonlinear function $\Theta$ is defined on the set of $G_2$-structures and has values in $\Omega^4(X)$, with  
\begin{equation}\label{theta}
\Theta(\psi) = {{*}_\psi}\psi,
\end{equation}
 where the Hodge star  ${{*}_\psi}$ is computed for the metric associated to the $G_2$-structure $\psi$.  
Proposition 10.3.5 in  \cite{Joyce} shows that if $\phi$ is a 
$G_2$-structure,  the derivative of $\Theta$ at $\phi$ satisfies
\begin{equation}\label{derivative}
\forall \psi \in \bz^3(M), d_\phi \Theta(\psi)  = {*}\mathcal P(\psi), 
\end{equation}
where the Hodge star corresponds to  $g$
and 
\begin{equation}\label{pe}
\mathcal P  =  \frac{4}{3} \pi_1
+\pi_7 - \pi_{27} .
\end{equation}
 Here $\pi_1$, $\pi_7$ and $\pi_{27}$ are the orthogonal
projections corresponding 
to the decomposition 
$\Lambda^3T^*M= \Lambda^3_1\oplus
\Lambda^3_7\oplus\Lambda^3_{27}$ associated to
the irreducible representations of $G_2$, see Lemma 3.2 in \cite{FeGr} or Proposition 10.1.4 in \cite{Joyce}.
Hence if $\psi = \dsurdt \phi(t)_{|t=0}\in \bz^3(M)$, we have
\begin{equation}\label{nabla}
 \nabla_{\psi} F = \sum_k (\eta_k\wedge \omega)\lrcorner 
{*} \mathcal P(\psi) \otimes \eta_k.
\end{equation}
\begin{lemma}\label{lemma surj}
The operator $\nabla F: \bz^3(M) \to \Gamma(Y,\nu)$ defined by equation (\ref{nabla})
is onto.
\end{lemma}
\bpr
Due to the properties of $\chi$, in this formula we can restrict 
our $\eta_k$'s to a local orthonormal frame of $\nu$ for the 
metric $g$. Now, recall \cite{FeGr} that $\Lambda^3_7= \{{*}(\phi\wedge \alpha), \alpha \in \Lambda^1T^*M\}.$
Consider $s\in \Gamma(Y,\nu)$, and $\alpha$
 the dual 1-form of $s$. More precisely, 
 $\alpha \in \Gamma(Y,T^*M)$ satisfies  
\begin{equation}\label{equasso}
\forall y\in Y, \ \forall v\in T_y M, \alpha_y (v) = \langle s(y),v\rangle.
\end{equation} 
 We choose $\omega$ such that $\phi (\omega) = 1$, which is always possible since $Y$ is associative. Since $\mathcal P$ acts as the identity on $\Lambda^3_7$  and ${*}$ is an involution, it is straightforward to see that 
 \begin{equation}\label{alpha}
\sum_l(\eta_l\wedge \omega) \lrcorner {*} \mathcal P({*}(\phi\wedge \alpha))\otimes \eta_l =  s.
\end{equation}
In order to prove the existence of $\psi \in \bz^3(M)$ such that
$\nabla_\psi  F = s$, we need to extend
${*}(\phi\wedge \alpha)$ outside $Y$ as a closed form. For this, let $p\in Y$, $U$ be an open set of $M$ containing $p$ and  local coordinates $y_1, y_2, y_3, x_1, x_2,x_3, x_4$ on $U$, where the $y_i$'s are 
coordinates on $Y$ and the $x_i$'s are
transverse coordinates. 
Because $Y$ is associative, the 3-form 
 $\psi' = {*}(\phi\wedge \alpha) \in \Gamma(Y,\Lambda^3 T^*M) $ is of the form
  $ \sum_{i=1}^4 dx_i\wedge \beta_i$ over $Y\cap U$, where for all $i$, $\beta_i$ is a 2-form. 
We  extend arbitrarily the $\beta_i$'s as  smooth 2-forms on $U$. 
Assume first that $s $ has  compact support in $U\cap Y$. Then so do the $\beta_i$'s on $U\cap Y$. Define 
$$\psi' = d(\chi_U\sum_i  x_i  \beta_i),$$
where $\chi_U$ is a cut-off function with support in $U$ and equal to $1$ in the 
neighbourhood of the support of $s$.  Then $ \psi'$ is a global closed 3-form with $\psi'_{|Y} = \psi$ and hence
satisfying $\nabla_{ \psi'} F = s$. For a general section $s\in \Gamma(Y,\nu)$, a partition of unity allows us to find $\psi \in \bz^3(M)$ 
such that $\nabla_\psi F = s$. We conclude that  $\nabla F$ 
is onto in the direction of $\bz^3(M)$. 
\epr
We can now finish the proof of Theorem \ref{generic}. If $ \bz_D$ is
the finite dimensional subspace of $\bz^3(M)$ generated by 
the former closed 3-forms $\psi$ associated to every $s\in \coker D$ given by Lemma \ref{lemma surj}, by the inverse 
mapping theorem, the set $$\bm = \{(\si, \psi)\in W^{k,p}(Y,\nu) \times \bz_D(M), F(\sigma, \psi) = 0\}$$ is a smooth manifold near $(0,\phi)$
if $k>1+3/p$. By the Sard-Smale theorem applied to the projection $\pi: \bm \to \bz_D$, for every generic $\psi\in \bz_D$ close enough to $\phi$, the slice $$\pi^{-1}(\psi)= \{\sigma \in W^{k,p}(Y,\nu), \exp_\sigma(Y) \text{ is $\psi$-associative}\}$$ is a smooth manifold or an empty set. As usual, the sections 
in $\pi^{-1}(\psi)$ are in fact smooth, hence the result.
\epr
\begin{remark}\label{torsion-free}  
By Theorem 10.4.4 in \cite{Joyce}, if $\phi$ is a torsion-free 
$G_2$-structure, the tangent space at $\phi$
of the set of torsion-free structures 
can be identified with $\mathcal L \oplus \mathcal H^3(M,\R)$, 
where $\mathcal L$ is the subspace of the Lie derivatives
of $\phi$, i.e. 
$\mathcal L =\{\mathcal L_X \phi, X\in C^0(M,TM)\}$, and $\mathcal H^3(M,\R) $ is 
the space of the real harmonic  3-forms on $M$. 
If $\psi =\mathcal L_X \phi \in \mathcal L$, Lemma  \ref{lemma Lie} below 
shows that 
the derivative of $F$ along $\psi$ equals $D X^\perp$, where $X^\perp \in \Gamma(Y,\nu)$
is the normal projection of $X$ onto the normal bundle of $Y$.  Hence, 
$\bl$ is of no use for $\nabla F$ to be onto. But the dimension of $\coker D$ is not in general 
less than $b^3(M)$, and even when it is, $\mathcal H^3(M) \to \coker D$
might well  be non injective (see the end of the subsection \ref{extension} for examples of every situation). 
This is the reason why we use the wider space of closed 
$G_2$-structures.
\end{remark}
\begin{lemma}\label{lemma Lie}
Let $M$ be a manifold equipped with a torsion-free $G_2$-structure $\phi$,  $Y$ be a smooth  compact closed $\phi$-associative submanifold and
$X$ a smooth vector field of $TM$ in the neighbourhood of $Y$. Then 
$$dF_{|(0,\phi)} (\mathcal L_X \phi) = D X^\perp,$$
where  $dF_{|(0,\phi)} (\mathcal L_X \phi)$ denotes
the derivative of the section $F$ given by (\ref{F psi}) at $(0,\phi)$ in the direction $\mathcal L_X \phi$,
$D$ is the Dirac-like operator given by (\ref{Operator}) and $X^\perp$ 
is the orthogonal projection of $X$ onto the normal bundle $\nu$ over $Y$. 
\end{lemma}
\bpr
Denote by $(\Phi_X^t)_{t\in [0,\epsilon]} $ the flow generated by $X$ near $Y$ and $\phi_t= \Phi_X^{t*} \phi$ the pull-back of $\phi$.
Hence, the metric $g_t$ associated to $\phi_t$ is $\Phi_X^{t*} g$, so that
$\Theta (\phi_t)= \Phi_X^{t*} \big( \Theta (\phi )\big),$
where $\Theta$ is defined by (\ref{theta}).
Let $(\eta_k^t)_{k=1, \cdots, 4}$ be an orthonormal framing of the normal bundle of $Y$ for the metric $g_t$, depending smoothly 
on $t$. Then
\beq
  F(0, \phi_t) &=&  \sum_k (\eta^t_k\wedge \omega)\lrcorner \Theta(\phi_t) \otimes \eta^t_k
 =   \sum_k (\eta^t_k\wedge \omega)\lrcorner \Phi^{t*}_X \big(\Theta(\phi)\big) \otimes \eta^t_k,
\eeq
 which implies
$$ dF_{|(0,\phi)} (\mathcal L_X \phi)  =  \sum_k (\eta_k^0\wedge \omega)\lrcorner \mathcal L_X \big(\Theta(\phi)\big) \otimes \eta_k^0$$
(there is no derivative of $\eta^t_k$ because $\omega \lrcorner  \Theta(\phi)=0$).  This is 
the right hand side of (\ref{derivee}) with $X$ instead of $s$. The end of the proof of Theorem \ref{Dirac} 
shows that $ dF_{|(0,\phi)} (\mathcal L_X \phi)  = D X^\perp$. 
\epr

In the following Proposition \ref{nonvide}, we give a situation where
we can find a way to isolate an associative after perturbing 
the $G_2$-structure.
\begin{prp}\label{nonvide}
Let $Y$ be a smooth closed $\phi$-associative submanifold, such that $\ker D$
is generated by a non vanishing normal vector field. Then there is 
a neighbourhood $V$ of $Y$ and 
a closed perturbation $\psi$ of $\phi$, such that the only element of $\bm_{Y,\psi}$ lying in $V$ is $Y$.
\end{prp}
\bpr
Let $\xi_1 \in \ker D \setminus \{0\}$. Since the normal bundle $\nu$ of $Y$ is trivial, we can find 
normal vector fields $\xi_2, \xi_3$ and $\xi_4$ such that the $\xi_i$'s form a global  framing of $\nu$. Let $(x_i)_{i=1, \cdots 4}$ 
the coordinates near $Y$ defined by exponentiating the $\xi_i$'s. The 3-form $\xi_1 \lrcorner {*}\phi$
writes $\sum_{i=2,3,4} dx_i\wedge \beta_i$. 
The closed
form $$\psi = d(x_1\sum_{i=2, 3,4} x_i\wedge \beta_i).$$
  is  defined near $Y$ and vanishes on $Y$. If $\phi_\lambda = \phi + \lambda \psi$, denote
by $g_\lambda$ the associated metric and by $D^\lambda$  the Dirac-like operator associated to $\phi_\lambda$.
We will prove that for $\lambda$ small enough, the only solution to $D^\lambda s = 0$ is the 
null section, which will prove Proposition \ref{nonvide} by Proposition \ref{smoothness}.

Derivative (\ref{derivee}) together with equation (\ref{derivative}) giving the derivative of 
the Hodge star imply that for every $s\in \Gamma(Y,\nu)$,
\begin{eqnarray*}
D^\lambda s &=&  
  \sum_k \eta_k\wedge \omega \lrcorner \left( \bl_s \big(\Theta (\phi) +\lambda {*}\mathcal P (\psi ) + O(\lambda^2)\big) \right)\otimes \eta_k \\
& = & Ds + \lambda \sum_k \eta_k\wedge \omega \lrcorner \left( \bl_s \big( {*}\mathcal P (\psi ) \big) \right)\otimes \eta_k + O(\lambda^2 s),
  \end{eqnarray*}
  where ${*}$ is the Hodge star associated to $\phi$ and $\bp$ is given by (\ref{pe})
  (note that $\eta_k$ is an orthonormal framing for every $\lambda$ since $g_\lambda = g$ on $Y$). 
  In particular,  if $s \in \ker D^\lambda,$
  \begin{equation}\label{estim2}
  Ds = O(\lambda s).
  \end{equation}
   Near $Y$, we have 
$
{*}\mathcal P (\psi ) =   x_1 \xi_1\lrcorner {*}\phi  + \sum_{i=2,3,4} x_i {*}\mathcal P (dx_1\wedge \beta_i ) + O(x^2),
$
so that on $Y$,
\begin{eqnarray*}
\mathcal L_s \big( {*}\mathcal P (\psi ) \big)=   s_1 \xi_1\lrcorner {*}\phi + \sum_{i=2,3,4} s_i {*}\mathcal P (dx_1\wedge \beta_i ),
\end{eqnarray*}
where $s =\sum_{i=1,2,3,4} s_i \xi_i$. This implies
\begin{equation}\label{stock}
D^\lambda s = Ds + \lambda s_1 \xi_1 + O(\lambda (s_i)_{i=2,3,4}) + O(\lambda^2 s).
 \end{equation}
Since $D$ is a self-adjoint elliptic operator, 
there is a constant $a$ (depending on $\lambda$ and $s$) such that  $s- a\xi_1 = O(Ds)$, 
see Corollary 5.7 in \cite{LaMi} for instance. 
Assume now that $s\in \ker D^\lambda$. 
Then, estimate (\ref{estim2})  implies that 
$s-a\xi_1  = O(\lambda s),$
so that projecting onto the directions $\xi_1$ and $\xi_i$, $i=2,3,4$, 
we get 
\begin{eqnarray}
s_1-a &=& O(\lambda s) \text{ and } \\ 
s_i &= & O(\lambda s) \text{ for } i=2,3,4. \label{estimation si}
\end{eqnarray}
The first estimate gives ${\nabla} s_1 = O(\lambda s)$, and the second one together with
equation (\ref{stock}) implies
\begin{equation}\label{Ds2}
 Ds =- \lambda s_1 \xi_1 + O(\lambda^2 s).
 \end{equation}
This gives 
\begin{equation}\label{D2s}
D^2s = - \lambda {\nabla} s_1 \times \xi_1 + O( \lambda^2 s) = O(  \lambda^2 s).
\end{equation}
Since $\ker D^2  = \ker D$ and $D^2$ is elliptic, we have by Corollary 5.7 in \cite{LaMi} and relation (\ref{D2s}) 
the estimate
$s-a\xi_1 = O(D^2 s )= O( \lambda^2 s) $, so that 
\begin{eqnarray}\label{last one}
Ds = O( \lambda^2 s) 
\end{eqnarray} 
 since $D\xi_1 = 0.$
Now, from (\ref{Ds2}) and (\ref{last one}) we deduce
$\lambda s_1 \xi_1 = O(\lambda^2 s)$. Since by (\ref {estimation si}), the norm of $s$ is equivalent to the norm of $s_1$ 
 when $\lambda $ tends to zero, the last estimate is impossible for
$\lambda $ small enough and a non zero $s\in \Gamma(Y,\nu),$ so that $s= 0$.  
\epr
This situation arises in particular in the Calabi-Yau extension, see Corollary \ref{coro pro S1 close} below.

\subsection{A vanishing theorem}
We turn now to the second way of getting the smoothness of the moduli space, namely Bochner's technique and Simons's theorem. We  formulate the following theorem which can be deduced 
from Theorem \ref{Simons}, since any associative submanifold is minimal.
\begin{thm}\label{closed}
Let $Y$ be a smooth closed compact  associative submanifold of a manifold $M$ with a closed 
$G_2$-structure. If the spectrum of $\br_\nu = \br - \mathcal A $ is positive, then 
$Y$ is isolated as an associative submanifold.
\end{thm}
 For the reader's convenience, we give below a proof of this result in the case where the $G_2$-structure is torsion-free. We will compute $D^2$ to use Bochner's technique. 
For this, we introduce the normal equivalent of the invariant second derivative. More
precisely, for any local vector fields $v$ and $w$ in $\Gamma(Y,TY)$, let $\nabla^{\perp 2}_{v,w}$ be the operator defined by
 $\nabla^{\perp 2}_{v,w} = \nap_v\nap_w - \nap_{\nabla^\top_v w}$
acting on $\Gamma(Y,\nu)$. It is straightforward to see that it is tensorial in  $v$ and $w$.
 Moreover, define the equivalent of the connection Laplacian: 
$$ \nabla^{\perp *}\nap = - \text{trace } (\nabla^{\perp 2}) = - \sum_i \nabla^{\perp 2}_{e_i,e_i},$$ 
where the $e_i$'s define a local orthonormal frame of $TY$.
\begin{thm}\label{D^2}
For $Y$ an associative submanifold in a manifold with a torsion-free $G_2$-structure, 
$ D^2 =  \nabla^{\perp *}\nap + \br_\nu.$
\end{thm}
We refer to the appendix for the proof of this theorem.
\bpr[ of Theorem \ref{closed}]
Let assume  that we are given a fixed closed associative submanifold $Y$. 
Consider a section $s \in \Gamma(Y,\nu)$. By classical computations using normal coordinates and thanks to Theorem \ref{D^2},
we have 
$$-\frac{1}{2}\Delta |s|^2  =  \sum_i\langle \nabla_i^\perp s, \nabla_i^\perp s\rangle   + \langle s, \nabla_i^\perp \nabla_i^\perp s\rangle  
                                =  |\nabla^\perp s|^2 - \langle D^2s, s\rangle   + \langle \br_\nu s, s\rangle . $$
 Since the Laplacian  equals  $-\text{div}  (\grad)$, its integral over the closed $Y$ 
 vanishes. We get:
\begin{eqnarray}\label{integral}
0 = \int_Y |\nap s|^2 - \langle D^2 s, s \rangle   + \langle \br_\nu s, s\rangle   dy.
\end{eqnarray}
Assume that    $s$ belongs to $\ker D$. Under the hypothesis that
 $\br_{\nu}$ is positive, the last equation implies  $s =0$. 
 Hence $ \dim \coker D = \dim \ker D = 0$, and by Proposition \ref{smoothness}, 
 $\bm_Y$ is a smooth manifold near $Y$ with vanishing dimension. In particular, $Y$ is isolated.  
\epr 

\section{Associative submanifolds with boundary}
In this section we explain our results in  
the case of an associative submanifold with boundary in 
a coassociative submanifold.  We first give below the principal results of \cite{GaWi}. For this, recall  that 
in a manifold with a $G_2$-structure and an associated  vector product $\times$, given $x\in M$ and $n$  an unit vector in $ T_xM$,  the application 
$$n\times: T_xM \to T_xM, v \mapsto n\times v$$ defines a complex structure on $n^\perp$, the orthogonal complement of $n$.  A 2-plane $L \subset n^\perp$ invariant under $n\times $ will be called a $n\times$-\textit{complex line}. 
\begin{thm}[\cite{GaWi}]\label{boundary} Let $M$ be a manifold
equipped with a $G_2$-structure $(\phi,g)$ and $Y$ a smooth compact associative submanifold with boundary in a coassociative submanifold $X$. 
Let $\nu_X$ be the normal complement of $T\partial Y$ in $TX_{|\partial Y}$, and 
$n$ the inward unit  normal  vector to $\partial Y$ in $Y$. Then
\begin{enumerate}
\item \label{un}
the bundle $\nu_{X}$ is a subbundle of  $\nu_{|\partial Y}$ and is a $n\times$-complex line, as is
the orthogonal complement $\mu_X$ of $\nu_X$  in $\nu_{|\partial Y}$. 
\item Viewing  $T\partial Y$, $\nu_X$ and $\mu_X$ 
as $n\times$--complex line bundles, we have $\mu_X^*\cong\nu_X\otimes_{\C}T\partial Y$.
\item Further, the problem of the associative deformations of $Y$ with boundary in $X$ is elliptic and of index
$\ind (Y,X) = \ind \overline \partial_{\nu_X} = c_1(\nu_X) + 1- g,$
where $g$ is the genus of $\partial Y$.
\end{enumerate}
\end{thm}
\begin{prp} \label{B-smoothness} 
Let $M$ be a smooth manifold equipped with 
a $G_2$-structure $(\phi,g)$ and let $Y$ be a smooth compact associative submanifold with boundary
in a coassociative submanifold $X$. Consider the adapted version of the linearization of (\ref{Operator}) for our boundary problem: 
$$ D: \beee_X=\{s\in \Gamma(Y,\nu), s_{|\partial Y} \in \nu_X\} \to \Gamma(Y,\nu).$$
If the cokernel of $D: \beee_X \to \Gamma(Y,\nu)$ vanishes, 
then $\bm_{Y,X}$ is smooth near $Y$ and of dimension equal to $\ind (Y,X)$. 
\end{prp}
\bpr
 For $2k>  3$ and $(k-r)/3 >   1/2$, define the adapted Banach space $\beee_X$ by  $$\beee_X = \{\sigma \in W^{k,2} (Y,\nu), \forall y\in \partial Y, \sigma(y)  \in \nu_{X,y} \}$$
and $ \bfff$ the bundle over $\beee_X$, where the fiber $\bfff_\sigma$ denotes $W^{k-1,2}(Y, \nu_\sigma)$. As before  $ \nu_\sigma$ is the normal bundle
to $ \exp_\sigma (Y)$. Let us assume first that $X$ is totally geodesic for the metric $g$. Then $\beee_X$ parametrizes the submanifolds
with boundary in $X$ and close enough to $Y$. 
Define the analogue of the map  (\ref{F}) in the proof of Theorem \ref{Dirac} by
$F: \beee_X \to \bfff,$  
$F(\sigma) =  \exp_\sigma^*\chi $. By the proof of 
Theorem \ref{Dirac}, $F$ is smooth and its derivative at the vanishing section
 is
$ D:\beee_X \to \Gamma (Y,\nu).$
Further, by Theorem 20.8 of \cite{BoWo}, 
the operator $D:\beee_X \to \Gamma (Y,\nu)
$ is Fredholm and Theorem  \ref{boundary} gives
its index. Now, if the cokernel of $D$ vanishes, 
then the inverse mapping theorem shows
that $\bm_{Y,X}$ is smooth near $Y$ and of 
the expected dimension equal to $\ind (Y,X)$. Lastly, Theorem 19.1 in \cite{BoWo} shows
that in fact, the sections belonging to $\bm_{Y,X}$
are smooth and so are the associated deformations of $Y$. 
In general, $X$ is not totally geodesic and as explained
in \cite{Bu} and \cite{KoLo}, $\exp_\si(\partial Y)$ 
has no reason to lie in $X$. For this, we 
change  the metric near $X$, as in the mentioned works. 
\begin{lemma}\label{hat} There exists  a tubular neighbourhood $U$ of $X$ and  a metric $\hat g$ such that
$\hat g(x) =g(x)$ for every $x\in X$,   $\hat g$ equals $g$ outside $U$,  and  $X$ is totally geodesic for $\hat g$. 
\end{lemma}
%
\bpr 
The exponential gives a diffeomorphism $\Phi$ between 
a tubular neighbourhood $U$  of $X$ in $M$ and a neighbourhood $V$ of the vanishing section in the normal vector bundle $N_X$
of $X$. Moreover, it sends $X$ to the vanishing section. 
Consider on $V$ the metric $h  = \pi^*g_{|TX} \oplus g_N,$
where $g_N$ is the natural flat metric on the fibers induced by the metric $g$, $ g_{|TX}$ is the induced metric on $X$ and 
$\pi: N_X \to X$ denotes the natural projection.
Now $H = \Phi^*h$ is a metric on $U$, for which
$X$ is clearly totally geodesic.  Take $\chi$ a cut-off function with support in $U$, equal to 1 in a neighbourhood of $X$. 
Then $\hat g = \chi H + (1-\chi) g$ satisfies
all the  conditions of the lemma.  
\epr 
Consider $\hat \nu$ the normal bundle over $Y$ for the new metric $\hat g$. For every section $\sigma \in \Gamma(Y,\hat \nu)$ we  use the adapted function 
$\hat F(\sigma) = \widehat{ \exp_\sigma}^*\chi(\omega)$, 
where $\omega$ can be chosen as before and
$\chi$ is the form associated to $\phi$, but $\widehat{ \exp}$ is the exponential map for the new metric $\hat g$. The proof of 
Theorem \ref{Dirac}
shows that differentiating $\hat F$ in the direction of $s \in \Gamma(Y,\hat \nu)$  gives the same  result $\nabla_s \hat F = Ds \in \Gamma(Y,\nu)$, even
if $s$ does not belong to $\Gamma(Y,\nu)$. 
Now, given a bundle isomorphism between $\hat \nu $ and $\nu$, 
it is straightforward to see that the kernel and the cokernel of $\hat \nabla F$ are isomorphic to the ones of $D$. The former conclusion in the totally geodesic case  still holds.
\epr

 \subsection{Varying the coassociative submanifold}
In subsection \ref{varying G2}, we perturbed the
$G_2$-structure in order for the moduli 
space $\bm_Y$ to become smooth. When the associative
submanifold has a boundary, we 
can repeat the same arguments. 
We can also move the boundary condition. 
As explained in the introduction, we will perturb generically $X$ as 
a smooth $\phi$-free submanifold, and no longer as a coassociative one.\\

\textbf{Theorem \ref{coasso} } 
\textit{Let $Y$ be a smooth associative submanifold
with boundary in a  smooth coassociative submanifold 
$X$. If the virtual dimension of $\bm_{Y,X}$ is non-negative, then for any sufficiently small generic smooth
deformation $X'$ of $X$,  there exists a small associative deformation $Y'$ of $Y$ such that $\bm_{Y',X'}$  
is smooth near $Y'$ and of  dimension equal to 
the index computed for the unperturbed situation. }
\bpr 
Recall \cite{McL} that
if $X$ is a coassociative submanifold, then its normal bundle
$N_X$ can be identified with the space of its self-dual two-forms
$\Omega^2_+(X)$. 
For $\alpha \in \Omega^2_+(X)$, define $\sigma_\alpha\in \Gamma (\partial Y, N_X)$ the restriction to $\partial Y$ of the associated normal vector field along $X$. By Theorem \ref{boundary},  $N_{X|\partial Y} = n \R \oplus \mu_X$, with
$n$ the inward unit normal vector to $T\partial Y$ in $TY$.   
Consider the subspace  $$\bc = \{\alpha \in \Omega_+^2(X), \sigma_\alpha \in \Gamma (\partial Y, \mu_X)\}.$$
Note that infinitesimal deformations of $X$ in these directions
are normal to $Y$. This will be considered
as the parameter space. 
For every $\alpha \in \bc$, extend $\sigma_\alpha$ 
to $\Gamma (Y,\nu)$ in the following way. The associative  $Y$ is diffeomorphic to  $Y_\epsilon = \partial Y \times [0,\epsilon]$ near
$\partial Y$, where $\partial Y$ holds for  $\partial Y \times \{0\}$. This allows us to identify 
$\nu_{|Y_\epsilon} $ with $\nu_{|\partial Y} \times [0,\epsilon]$
and so this gives an extension 
of  $\sigma_\alpha$ on $Y_\epsilon$. 
Take $\rho$ a cut-off function satisfying $\rho= 1$ in 
the neighbourhood of $\partial Y$ and 
with support in $Y_\epsilon$. 
Then $\hat \sigma_\alpha = \rho \sigma_\alpha\in \Gamma(Y,\nu) $ 
is a smooth normal vector field along $Y$ such that $\hat \sigma_\alpha = \sigma_\alpha $ near $\partial Y$. 
Now, let $\beee_\partial$ be the set 
$$ \beee_\partial = \{(\alpha, s)  \in \bc\times  \Gamma(Y,\nu), 
\forall y\in \partial Y, s(y) \in T_yX   \}.$$
Here we will assume that $X$ is totally geodesic 
as in the first part of the proof of Proposition \ref{B-smoothness}. If not, we change the metric by Lemma \ref{hat}. Hence if $(\alpha,s)\in \beee_\partial$ and if we define  $\phi_{\alpha, s} = \exp_{\hat \sigma_\alpha} \circ \exp_s$, then $Y_{\alpha, s} =  \phi_{\alpha, s} (Y) $ is a smooth submanifold with boundary in 
$X_\alpha = \exp_{\sigma_\alpha}(X)$. Let  $\bfff$ be the bundle over $\beee_\partial$, where the fiber $\bfff_{\alpha,s}$ equals $\Gamma(Y,\nu_{\alpha,s} )$ and $\nu_{\alpha,s} $ denotes
 the normal bundle of
$ Y_{\alpha,s}$. Define the section $F: \beee_\partial  \to \bfff $ by 
$F(\alpha,s)  =  \phi_{\alpha,s}^* \chi(\omega).$
Then $Y_{\alpha, s} $ is an associative
submanifold if and only if $F(\alpha,s) = 0$. 
Now for every fixed $\alpha \in \bc$, consider
the restriction map 
\beq 
F_\alpha: \{s \in \Gamma(Y,\nu), 
s_{|\partial Y} \in TX   \} & \to &\Gamma (Y,\nu_{\alpha,s})\\
s &\mapsto &F(\alpha,s)
\eeq
Two  tedious computations analogous to the proof of Theorem \ref{Dirac} and the proof of Theorem \ref{boundary} in Section 4 of \cite{GaWi}  show that for every $\alpha \in \bc$, the derivative of $F_\alpha$ is elliptic in the sense of  Definition 18.1 of \cite{BoWo}. Further, $F_\alpha$ is
clearly  a deformation of $F_0$, hence 
 $F_\alpha$ is a Fredholm map of index computed in Theorem \ref{boundary}. For a genericity result, we need the  classical 
\begin{thm}[Theorem 1.5.19 of \cite{Nico}]
Let $\bc$, $\mathcal E$ and $\mathcal F$ be Banach spaces, $F: \bc \times \mathcal E \to \mathcal F$
a smooth map, such that for every $\alpha \in \bc$, 
$F_\alpha = F(\alpha,.)$ is a Fredholm map between $\mathcal E$ and $\mathcal F$. 
If $dF: \bc \times \mathcal E \to \mathcal F$ is onto
at $(\alpha_0, x_0)$, then $F^{-1}(y_0)$ 
is locally a smooth manifold, where $y_0= F(\alpha_0,x_0)$. Further, for every generic $\alpha \in \bc$ close
enough to $\alpha_0$, the fiber $F_\alpha^{-1}(y_0)$ is a smooth manifold
of finite dimension equal to  the index of 
$F_\alpha$. 
\end{thm}
We compute 
the derivative of $F$ at $(0,0)\in \beee_\partial$. 
One can easily check using the proof of Theorem  \ref{Dirac} that this is equal to
\beq \nabla_{(0,0)} F:
\beee_\partial  & \to & \Gamma(Y,\nu)\\
(\alpha,s) & \mapsto & D(s+\hat \sigma_\alpha).
\eeq
This derivative is onto. Indeed,  let $s'$ be
a section in $\Gamma(Y,\nu)$. Since $Y$ has a boundary, our Dirac-like operator $D$ is onto by Theorem 9.1 of the 
book \cite{BoWo}, so there is a section $s\in \Gamma(Y,\nu)$ such that $Ds = s'$. Now decompose $s_{|\partial Y}$
as  $s_\nu +s_\mu $ with $s_\nu \in \Gamma(\partial Y, \nu_X)$
and $s_\mu \in \Gamma(\partial Y, \mu_X)$. 
Choosing the 2-form $\alpha \in 
\bc$ such that  $s_\mu = \sigma_\alpha$, we have 
$D((s-\hat \si_\alpha)+ \hat \sigma_\alpha) = s'$ with $(\alpha, s- \hat \si_\alpha)\in \beee_\partial$, hence the result.
\epr
 As in Theorem \ref{generic}, 
we can restrict our smoothing deformations to a finite dimensional space of dimension equal to $\dim \coker D$.

 \subsection{A vanishing theorem}
Given $Y$ an associative submanifold with boundary in a coassociative submanifold $X$, we turn now to metric conditions on $Y$  that insure local smoothness
 of the moduli space $\bm_{Y,X}$. Let $\nu$ be the normal bundle of $Y$ and $n$ is the inward normal vector to $\partial Y$
in $Y$.  Recall that if $L\subset \nu$ is a $n\times$-complex line bundle over $\partial Y$, the  operator  $\bd_{L}: \Gamma(\partial Y, L) \to \Gamma(\partial Y, L) $  was defined in the introduction by $\bd_L s =   \pi_L(v\times \nabla^\perp_w s- w\times \nabla^\perp_v s),$
where $\pi_L: \nu \to L$ is the orthogonal projection to $L$ and $\{v,w= n\times v\}$ a local orthonormal frame for $T\partial Y$. 
We refer to the appendix for  the proof of the following proposition.
\begin{prp}\label{bd}
The operator  $\bd_L$ is of order 0, symmetric, and its trace is $2H$, where $H$ is the mean curvature of $\partial Y$ in $Y$ with respect to $-n$.
\end{prp}
Moreover, consider the operator 
$ (D,L)$ defined by $ D: \{s\in \Gamma(Y,\nu), s_{|\partial Y} \in L\} \to \Gamma(Y,\nu) $. 
We will use the following lemma, whose proof can be found  in the appendix. 
\begin{lemma} \label{adjoint} We have 
$\coker (D,L) = \ker (D,L^\perp),$ where $L^\perp$ is
the orthogonal complement of $L$ in $\nu_{|\partial Y}$. 
\end{lemma}
We now prove the vanishing theorem stated previously:

\textbf{Theorem \ref{open} } \textit{Let $M$ be a manifold equipped with a torsion-free $G_2$-structure and  $Y$ be an associative submanifold  with boundary in a coassociative submanifold $X$. 
If  $\bd_{\mu_X}$ and 
$\br-\mathcal A $ are positive, the moduli space $\bm_{Y,X}$ is  smooth near $Y$ and of dimension 
given by the virtual one.}
\bpr 
To prove Theorem \ref{open}, it is enough by Proposition \ref{B-smoothness} to show that 
$ \coker (D,\nu_X)$, which equals $  \ker (D,\mu_X)$ by Lemma \ref{adjoint}, is trivial. So let $s \in   \ker (D,\mu_X).$   Since $Y$ has a boundary, we need to change the integration (\ref{integral}), because the divergence has to be considered: 
\begin{eqnarray}
 \int_Y |\nap s|^2  + \langle \br_\nu s, s\rangle   dy = \frac{1}{2} \int_Y \text{div }\nabla |s|^2 dy.
 \end{eqnarray}
By Stokes, the last  equals 
\beq
                             -\frac{1}{2}\int_{\partial Y} d|s|^2(n) d\sigma =  -\int_{\partial Y} \langle \nap_ns, s\rangle   d\sigma,
                                                                  \eeq 
where $n$ is the inward   unit normal vector of $\partial Y$. Choosing a local
orthonormal frame $\{v, w = n\times v\}$ of $T\partial Y$, 
$0= Ds = n\times \nabla^\perp_n s + v\times \nabla^\perp_v s+ w\times \nabla^\perp_w s$
implies that $$\nabla^\perp_n s =- w\times \nap_v s+ v\times \nap_w s.$$ Here we used the formula 
\begin{eqnarray*}
 \forall u,v,w \in TM, \ \chi (u,v,w) = -u\times (v\times w) - \langle u,v\rangle w + \langle u,w\rangle v,
 \end{eqnarray*}
  so that  for orthogonal vectors $ u,v,w\in TM,$
 \begin{eqnarray}
 u\times (v\times w)&=& w\times (u\times v) \text{ and} \label{permut} \\
u\times (u\times w)&= &-||u||^2 w. \label{permut1}
 \end{eqnarray}
Hence, 
\begin{eqnarray*}
     -\int_{\partial Y} \langle \nap_ns, s\rangle   d\sigma =            \int_{\partial Y} \langle w\times \nap_v s- v\times \nap_ws, s\rangle   d\sigma
           = -\int_{\partial Y} \langle \bd_{\mu_X} s,s\rangle   d\sigma.
\eeq
Summing up, we obtain the equation 
\begin{equation}\label{integration}
\int_Y |\nap s|^2 dy + \int_Y\langle \br_\nu s, s\rangle   dy + \int_{\partial Y} \langle \bd_{\mu_X} s,s\rangle   d\sigma=0.
\end{equation}
If $\bd_{\mu_X}$ and $\br_\nu$ are
 positive,  $s$ vanishes, hence the result.
\epr
\section{Examples}\label{exemples}
\subsection{Flatland}
In flat spaces, the curvature tensor $R$ vanishes, and so  $\br_\nu = -\ba\leq 0$. Consequently, a priori Theorem \ref{open} does not apply. Nevertheless, we have the 
\begin{cor}\label{flat}
Let $M$ be a manifold equipped with a torsion-free $G_2$-structure whose metric is flat, and $Y$ be a totally geodesic associative submanifold with boundary
in a coassociative $X$. If  $\bd_{\mu_X}$ is
 positive, then $\bm_{Y,X}$ is  smooth near $Y$ and of the expected dimension.
\end{cor}
\bpr The hypotheses on $M$ and $Y$ imply that $\br_\nu =0$.
Consider $s \in \coker (D,\nu_X) = \ker (D,\mu_X)$. 
Formula (\ref{integration}) shows 
that $\nabla^\perp s = 0$ and $s_{|\partial Y} = 0$. 
Using $d|s|^2 = 2\langle \nap s , s\rangle   = 0$. This implies $s = 0$ and the result. 
\epr
When $M = \R^7$ with its canonical flat metric, we get the following very explicit example considered in \cite{GaWi}. 
 Take a ball $Y$ in $\R^3\times \{0\}\subset \R^7$ with real analytic boundary, 
and choose  any normal real analytic vector field  $e\in \Gamma(\partial Y, \nu )$.  By \cite{HaLa}, there is a unique local coassociative $X_e$ containing $\partial Y$  such 
that its tangent bundle $T_yX_e$ contains $e(y)$ at every boundary point $y$.
\begin{moncor}\label{ball}
Let us assume that  $Y$ is a strictly convex ball in $\R^3$. Then there exists a positive constant $\epsilon$, such that for every normal vector field $e\in \Gamma(\partial Y, \nu )$ satisfying  $||de||_{L^\infty} \leq \epsilon$, the moduli space $\bm_{Y,X_e}$ is smooth near $Y$ and 
one dimensional.   
\end{moncor}
\bpr Since the fiber bundle $\nu_{X_e}$ is trivial and the genus of $\partial Y$ is zero, the index equals here $c_1(\nu_X)+1-g = 1$.
 We want to show that $\bd_{\mu_X}$ is positive. To see that, 
 we choose local orthogonal principal directions $v$ and $w = n\times v$ 
 in $T\partial Y$. From Theorem \ref{boundary}, 
 we know that  $v\times e$ is a non vanishing section of $\mu_X$. 
Let assume first that $e$ is constant. We compute, using relation (\ref{permut}),
\beq
\bd_{\mu_X}(v\times e) &  = & v \times ( \nabla^{\perp \partial}_w v \times e) - 
                                                     w  \times (\nabla^{\perp \partial}_v v \times e)\\
                                         &  = & - k_v w\times (n\times e) =  k_v v \times e,
\eeq
where $k_v$ is the principal curvature in the direction of $v$.
This shows that $k_v$ is an eigenvalue of $\bd_{\mu_X}$, and 
since we know that its trace is $2H$ by  Proposition \ref{bd}, we get that the other eigenvalue is $k_w$, the other principal curvature of $\partial Y$. 
These eigenvalues are positive if the boundary of $Y$ is strictly convex and Corollary \ref{flat}
gives the result.Now, if
$e$ is close enough to be a constant vector field, the eigenvalues of $\mathcal D_{\mu_X}$ remain positive, hence the general result.
\epr
 In fact, in the case where $e$ is constant, we can give a better statement.  Indeed,
let $s \in  \ker (D,\nu_X)$, 
and decompose  $s_{|\partial Y}$ as $s = s_1 e + s_2 n\times e$.  
Of course, $e$ is in the kernel of $\bd_{\nu_X}$, and hence 
by Proposition \ref{bd}, the second term is an eigenvector of $\bd_{\nu_X}$ for the  eigenvalue $2H$. So
 formula (\ref{integration}) applied to $s$ gives
$ \int_{Y} |\nap s|^2+ \int_{\partial Y} 2H|s_2|^2 = 0.$
If $H>  0$, this implies immediately that $s_2=0$ and $s_1$
is constant, so $s$ is proportional to $e$. This proves that
$\dim \ker (D,\nu_X) = 1$ under the weaker condition that $H>  0$. Lastly, in fact we can even show that  
$\bm_{Y,X_e}=\R$.

\subsection{The Bryant Salamon construction}
\textbf{The spin bundle and its metric. }As recalled briefly in the introduction, Bryant and Salamon
\cite{BrSa} found on the total spin bundle $\bs \simeq S^3 \times \R^4$ of the 
round sphere $S^3$ a complete metric with holonomy precisely equal to $G_2$.
This metric is of the form $$g =\alpha(r) \pi^* g_S + \beta(r) g_v,$$ 
where $g_v$ is the
flat metric on the fiber $\bs_x \simeq \R^4$ induced by $g_S$, $r$ is its associated norm, $g_S$ the round metric on $S^3$ and $\pi: \bs \to S^3$ the natural projection. For some particular smooth functions $\alpha$ and $\beta$, the authors
proved that the holonomy of the metric is $G_2$. 
In this situation, the base $S^3$ is associative and the Dirac operator 
$D$ is the classical one for the spin bundle $\bs$. 
\begin{cor}[\cite{McL}]
The associative $S^3$ is isolated as an associative submanifold. 
\end{cor}
\bpr By the famous computation of Lichnerowicz \cite{Li}, 
$D^2 = \nabla^*\nabla  + s/4,$
 where $s$ is the scalar curvature of $(S^3,g_S)$ and $\nabla $ is the induced connection on the spin bundle, which is in our case the connection $\nabla^\perp$. Identifying with the equation 
in Theorem \ref{D^2}, we get that $\br_\nu = s/4$. Since $S$ is positive, so is $\br_\nu$, and Theorem \ref{Simons} then implies the result.  \epr

\noindent
\textbf{Example with boundary.} Choose a point $p$ on the base $S^3$, 
a ball $B_\rho \subset \bs $ of  radius  $\rho$ around $p$ and define 
$ Y_\rho = B_\rho\cap S^3$. Take  a normal vector field $e\in \Gamma (\partial Y_\rho, \nu)$ at the boundary of the associative $Y_\rho$. Here $\nu_y = \bs_y$ for $y\in \partial Y_\rho$. The round sphere is real algebraic as is its metric $g_S$, hence we can find for $\rho $ small enough a local chart $\Phi: B_\rho \to \R^7$ such that
 $\Phi(Y_\rho)\subset \R^3\times \{0\}$, and $\Phi_*g$ is a real analytic metric. Further we choose $B_\rho$ and $e$ in such a way that  $\Phi(\partial Y_\rho)$ and $\Phi_*e $ are real analytic. Now, a straightforward generalization
of the arguments in \cite{HaLa} based on the Cartan-K\"{a}hler 
theory proves that $e$ and $\partial Y_\rho$ generate a semi local coassociative submanifold $X_e$  containing $\partial Y_\rho$. 
\begin{moncor}\label{BrySal}
For $\rho $ small enough, $\bm_{Y_\rho, X_e}$
is  smooth near $Y_\rho$  and one dimensional. 
\end {moncor}
\bpr The genus of $\partial Y_\rho$ vanishes and the 
subbundle $\nu_{X_e}$ is trivial, hence the index of the associative deformations problem equals one. 
We
can assume that $\Phi_*g(0) $ is the standard metric of $\R^7$, hence   $d_p\Phi(\bs_p) = 0\oplus \R^4$. Moreover we choose $\Phi$ such that the Levi-Civita connection
of $\Phi_* g$  vanishes at $0$. When $\rho$ tends to zero, $\Phi(\partial Y_\rho)$ is asymptotically close to be the round ball $\rho B^3 \subset \R^3$ for the metric $g_0$.  Then we know from the proof of Corollary \ref{ball} that the eigenvalues of the operator $\bd_{\mu_{X_e}}$ computed in  the model situation (i.e. with the flat metric and connection) equal the principal curvatures, here the inverse of $\rho$. Hence for $\rho $ small enough, $\bd_{\mu_{X_e}}$ and $\br_\nu = s/4$ are both positive.  Theorem \ref{open} then implies the result.
\epr

\subsection{The Joyce construction}
Recall briefly the construction of the compact smooth manifold with holonomy $G_2$ constructed by Joyce in section 12.2  of \cite{Joyce} 
and used in \cite{GaWi} for an example of an associative with boundary. 
 On the flat torus $(T^7,g_0)$ equipped with the $G_2$ structure 
$\phi_0 = dx_{123}+ dx_{145}+ dx_{167}+ dx_{246}- dx_{257} - dx_{347}-dx_{356},$
let 
\beq
\alpha : (x_1, \cdots, x_7) \mapsto (x_1, x_2, x_3, -x_4, -x_5, -x_6, -x_7),\\
\beta :  (x_1, \cdots, x_7) \mapsto (x_1, -x_2, -x_3, x_4, x_5, \frac{1}{2}-x_6, -x_7),\\
\gamma :  (x_1, \cdots, x_7) \mapsto (-x_1, x_2, -x_3, x_4, \half -x_5, x_6, \half -x_7),\\
\sigma_0: (x_1, \cdots, x_7) \mapsto (x_1, \half -x_2, \half -x_3, x_4, x_5, -x_6, \half -x_7),\\
\tau_0 : (x_1, \cdots, x_7) \mapsto (x_1, x_2, \half -x_3, \half - x_4, x_5, x_6, \half -x_7)
\eeq
be isometric involutions,  where $\si_0^*\phi_0 = \phi_0$ and $\tau_0^*\phi_0 = -\phi_0$. If $\pi : T^7\to T^7/\Gamma$ is the quotient of $T^7$ by $\Gamma$  the group
generated by $\alpha$, $\beta$ and $\gamma$, one can check that the image $Y$ by $\pi$ of 
$\{(x_1, \quart, \quart, x_4, x_5, 0, \quart), x_{1,4,5}\in T^3\}$
in $T^7/\Gamma$ is a smooth closed associative submanifold $Y$ belonging to the fixed point set of the well defined
involution $\si =\pi_*\si_0$.  Likewise, the image $Y_\partial$ by $\pi$ of 
$\{(x_1, \quart, \quart, x_4, x_5, 0, \quart), x_{1,5}\in T^2, \quart \leq x_4\leq \troisquart\}$
is a smooth associative submanifold with boundary. This boundary is the union of two 2-tori 
 embedded in the two disjoint smooth coassociatives $X_1$ and $X_2$, where $X_i$ is the image by $\pi$ of
$\{(x_1, x_2, \quart, a_i, x_5, x_6, \quart), x_{1,2,5,6}\in T^4\}$ with $a_1 = \quart$ and
 $a_2 = \troisquart$. The latter submanifolds  are components
of the fixed point set of $\tau = \pi_*\tau_0$. Joyce's method to construct a metric with holonomy
precisely equal to $G_2$ on a resolution $M$ of the singularities of $T^7/\Gamma$ can be made $\si$- and $\tau$-equivariantly, so that
after the process $Y$, $Y_\partial$, $X_1$ and $X_2$  remain associative and coassociative, respectively. 
Now, the bundles $\nu_{X_i}$, $i=1,2$,  are clearly trivial over the two components of $\partial Y_\partial$, so that the index of the deformation problem vanishes. 
From Theorem \ref{generic} and Theorem \ref{coasso} we get that 
for every generic closed perturbation 
$\psi$ of the $G_2$-structure, $Y$ disappears or is perturbed into an isolated closed $\psi$-associative torus. 
 Likewise, for every generic small  $\phi$-free deformation $\tilde X_i$ of $X_i$ there is a perturbation $\tilde Y_\partial$ of $Y_\partial$ such that $\bm_{\tilde Y_\partial, \tilde X} $ is a singleton near $\tilde Y$ or is empty. 
 \begin{remark} We would like to know which  alternative holds. 
Unfortunately,  even if we are far from the singularities of  $T^7/\Gamma$,
we don't know how to improve Proposition 11.8.1 of \cite{Joyce} 
in order to get a control of the $\eta_j$'s in $C^2$-norm. Said otherwise,
 the perturbation of the metric 
 a priori has  effects on the whole $M$, and can be big in $C^2$  norm. Hence, our methods do not allow to understand 
 the effects of the perturbation on the associative submanifolds. 
\end{remark}


\subsection{ Extensions from the Calabi-Yau world}\label{extension}
 \textbf{The closed case. } Let $(N, J, \Omega, \omega)$  be a Calabi-Yau 6-dimensional manifold, 
where $J$ is an integrable complex structure, $\Omega$ a non vanishing holomorphic 3-form
and $\omega $ a \kah form. 
Then $M = N\times S^1$  
is a manifold with holonomy in $SU(3) \subset G_2$. An associated torsion-free $G_2$-structure on $M$  is given by  $ \phi = \omega \wedge dt + \text{Re } \Omega. $
Recall that a closed special Lagrangian $L$ in $N$ is a 3-dimensional submanifold  satisfying both 
conditions $\omega_{|TL} = 0$ and $\text{Im } \Omega_{|TL} = 0$.  We know from \cite{McL} that $\bm_L$ the moduli space of  special Lagrangian deformations of $L$ is smooth
and of dimension $b^1(L)$. Now for every $t\in S^1$, the product   $Y= L\times \{t\}$ 
of a special Lagrangian and a point is a $\phi$-associative submanifold of $M$. The following is inspired by a analogous result on coassociative submanifolds of Leung ($\cite{Le}$, Proposition 5): 
\begin{prp}\label{Calabi}
Let $t\in S^1$. The moduli space $\bm_{L\times \{t\}}$ of associative deformations of $L\times \{t\} $ is always smooth, and can be identified with the product $\bm_L\times S^1$, hence  of dimension $b^1(L) +1$.
\end{prp}
\bpr   
Consider a closed associative submanifold $Y$ in the same homology class as $L\times \{t\}$. On the one hand, $Y$ has a bigger volume than its projection $\pi(Y)$ to $N\times \{t\}$ and equality holds only if $Y$ lies in $N\times \{t'\}$ for a constant $t'$. 
On the other hand, $\pi(Y)$ is in the same homology class as $L$, hence has volume larger than that of $L$, since
special Lagrangians minimize the volume in their homology class. But $Y$ is associative, hence has the same volume as $L$. Consequently all these volumes equal, and $Y$ is of the form $L'\times \{t'\}$. It is now immediate that $\phi$-associativity of $Y$ implies that $L'$ is special Lagrangian. 
\epr
For the sequel, we will need another 
\bpr[ of Proposition \ref{Calabi}] 
 Recall that since $L$ is Lagrangian, its normal bundle $NL$ is simply $JTL$, and  the normal bundle $\nu$ of $Y =L\times \{t\}$ is isomorphic to $JTL\times \R \partial_t$,
 where $\partial_t$ is the dual vector field of $dt$. 
In this situation, we don't use the expression for  $D^2$ given in Theorem \ref{D^2}.
Instead, we give another formula for it. 
 If $s= J\sigma\oplus \tau \partial_t$ is a section of $\nu$, with $\sigma \in \Gamma(L,TL)$
 and $\tau \in \Gamma(L,\R)= \Omega^0(L)$, 
 we call $\sigma^\vee \in \Omega^1(L,\R)$ the 1-form dual to $\sigma$, and we use
 the same symbol for its inverse. Moreover, we use the classical notation  ${*}: \Omega^k(L) \to \Omega^{3-k}(L)$
 for  the Hodge star. Lastly, we define: 
 \beq 
 D^\vee: \Omega^1(L)\times \Omega^0(L) &\longrightarrow &\Omega^1(L)\times \Omega^0(L)\\
 (\alpha, \tau)& \mapsto& ((-J\pi_L D(J\alpha^\vee, \tau))^\vee, \pi_t D(J\alpha^\vee, \tau)),
 \eeq
  where $\pi_L$ (resp. $\pi_t$) 
 is the orthogonal projection $\nu = NL \oplus \R$ to the first (resp. the second) component. This is just a
 way to use  forms on $L$ instead of normal ambient vector fields.
\begin{prp}\label{harmonic}
For every $(\alpha, \tau) \in \Omega^1(L)\times \Omega^0(L)$, 
\beq 
D^\vee (\alpha,\tau) &=& (-{*}d\alpha - d\tau, {*}d{*}\alpha) \\
(D^\vee)^2 (\alpha, \tau) &=& - \Delta (\alpha, \tau),
\eeq
where $\Delta  = d^*d + dd^*$ (note that it is $d^*d$ on $\tau$). 
\end{prp}
We refer to the appendix for the proof of this Proposition.
We see that for an infinitesimal associative  deformation of 
$L\times \{t\}$, then $\alpha$ and $\tau$ are harmonic 
over the compact $L$.  In particular,  $\tau$ is constant and $\alpha$
describes an infinitesimal  special Lagrangian deformation of $L$ (see \cite{McL}).  
In other words, the only way to displace $Y$ is 
to perturb $L$ as special Lagrangian  in $N$ or translate
it along the $S^1$-direction. Lastly, $\dim \coker D = \dim \ker D= b^1(L) + 1 $ and by an immediate refinement of Proposition \ref{smoothness} for cokernels  with constant dimension,
$ \bm_Y$ is smooth and of dimension $b^1(L) + 1$. 
\epr

\textbf{Symmetry breaking.} Although the moduli space is smooth, 
the deformation problem for $L\times \{.\}$ is always obstructed. Theorem \ref{generic} proves that any closed generic perturbation of the $G_2$-structure $\phi$ will make the $S^1$-symmetry disappear as well
as the $\bm_L$-family of associative submanifolds. We give here a family of examples of this phenomenon:
\begin{cor} \label{coro pro S1 close}
Let  $L$ be  a smooth closed special Lagrangian sphere in $N$, $t_0\in S^1$ and 
$Y= L\times \{t_0\} $ in $ N\times S^1$ 
equipped with the $G_2$-form 
$\phi = \text{Re } \Omega + \omega\wedge dt$ and $f: S^1 \to \R$ a smooth function vanishing transversally at a finite number of points in $S^1$. Then, 
there is a closed perturbation $\psi$ of $\phi$  such that 
the  connected  components  of  $L  \times  f^{-1}(0)$
   are  associative  with respect to $\psi$,  and  are  the  only  $\psi$-associatives  near $ \{ L  \times  \{t\}  :  t  \in  S^1  \}$.
\end{cor} 
\bpr 
Define $\tilde \psi = -f(t) {*} (\phi\wedge dt) = -f(t) \dsurdt \lrcorner {*} \phi = f(t)\text{Im } \Omega $ 
on $L\times S^1$ since ${*}\phi = \text{Im } \Omega\wedge dt + \frac{\omega^2}{2}$. 
We extend $\tilde \psi$ as a closed 3-form $\psi$ following the proof of Lemma \ref{lemma surj}:
since $L$ is special Lagrangian, $\text{Im } \Omega_{|L}\in \Gamma(L,\Lambda^3 T^*N)$ can locally  be written as
$\sum_{i=1}^3 dx_i\wedge \beta_i$, 
where $(x_i)_{i=1,2,3}$ are local normal coordinates in $N$ over $L$. If $(\chi_U)_U$ is a finite set of cut-off functions
in a neighbourhood of $L$, then the closed 3-form $\psi = d(f(t)\sum_{U,i} \chi_U x_i \beta_i)$
is well defined on $N\times S^1$  and
satisfies $\psi = f(t) \text{Im } \Omega + O(dist(., L\times S^1))$. 
We choose as a closed perturbation the 3-form $\phi_\lambda = \phi +\lambda \psi$. 
Now take $t_0\in S^1$, such that $f(t_0)=0$. If we choose coordinates on $S^1$ 
such that $t_0=0$, then there 
exists $a\not=0$ with $f(t) = at+O(t^2)$.  
Proposition \ref{nonvide} shows that for $\lambda$ small enough,
$L\times\{t_0\}$ is the only local $\psi-$associative. 
Now take $t_0$ such that $f(t_0)\not= 0$. The following lemma holds in a general situation: 
\begin{lemma}\label{lemma empty}
Let $Y$ be a compact smooth associative submanifold of $M$ equipped with a closed $G_2$-structure $\phi$, such that
near $Y$, $\bm_{Y,\phi}$ is one-dimensional and $\dim \ker D = 1$ at every element of $\bm_{Y,\phi}$.  
Let $\xi\in \Gamma(Y,N_Y)$ be a non trivial normal vector field in  $\ker D $ and
$\tilde \psi$ be the 3-form $\xi\lrcorner {*} \phi\in \Gamma(Y,\Lambda^3T^*M)$. 
If $\psi$ is any closed extension of $\psi$ in a neighbourhood of $Y$ and 
$\phi_\lambda = \phi + \lambda \psi$, then for $\lambda \not= 0$
small enough the moduli space 
$\bm_{Y,\phi_\lambda}$ near $Y$ is empty. 
\end{lemma}
\bpr By definition of $\phi_\lambda$ and Lemma \ref{lemma surj}, the derivative of $F(\lambda, s) = \exp_s^*\chi_{\phi_\lambda}(\omega)$ is 
of index 1 and surjective
at $(\lambda=0,s=0)$, so that the vanishing locus of $F$ is locally smooth, of dimension 1 and 
contains $\bm_{Y,\phi} $. These sets must be locally equal, hence the result.
\epr
We come back to the situation described in Proposition \ref{coro pro S1 close} . 
If $t$ is such that $f(t)\not= 0$, Lemma \ref{lemma empty} shows that
$\bm_{L\times \{t_0\},\phi_\lambda}$ is empty for $\lambda $ small enough. 
\epr

\textbf{Coclosed deformations.} If we prefer \textit{coclosed} deformations of the $G_2$-structure we get a more precise statement and a very short proof:
\begin{prp} \label{prop S1 coclosed}
Let  $L$ be  a smooth closed special Lagrangian sphere in $N$, $Y= L\times \{1\} \subset N\times S^1$ and $f: S^1 \to \R$ a smooth function vanishing  at a finite number of points in $S^1$. For every $\lambda \in \R$,  define
$\phi_\lambda = \text{Re } (e^{i\lambda f(t)} \Omega) + \omega\wedge dt$ a family of coclosed $G_2-$structures.  Then,  if $\lambda \not= 0$, 
$\bm_{Y,\phi_\lambda} = f^{-1}(0)$ near $L\times S^1$. 
\end{prp} 
Note that in particular, the  transversality condition for $f$ is no more needed.
\bpr 
The proof is almost the same as the proof of Proposition \ref{Calabi}. Take $Y$ a $\phi_\lambda$-associative submanifold of $N\times S^1$
in the same class of homology as $L\times \{1\}$. 
Since the metric associated to $\phi_\lambda$ is independent of $\lambda$, the arguments of Proposition \ref{Calabi} still hold, and $Y$ writes
$L'\times \{t'\}$ for some submanifold $L'\in N$ and $t\in S^1$. The latter $L'$ must be a special Lagrangian for $e^{i\lambda f(t)}\Omega$ since
$Y$ is $\phi_\lambda$-associative. Hence, $\text{Im } (e^{i\lambda f(t)} \Omega )$ vanishes on $TL'$. But $L'$ lies in the same class of homology as $L$,
so $\int_L \text{Im }(e^{i\lambda f(t)}\Omega )$ should vanish because $\Omega $ is closed. Now, this is in fact $\int_L  \sin(\lambda f(t))\text{Re }\Omega = 
\sin(\lambda f(t)) Vol (L)$ which is non zero if $\lambda \not= 0$ is small enough (independently of $t$) and $f(t)\not=0$. If $f(t) =0$ 
and $L'$ is close enough to $L$, then $L'=L$  since a special Lagrangian sphere is isolated. 
Note that $\phi_\lambda$ is coclosed because $*\phi_\lambda = \text{Im } (e^{i\lambda f(t)}\Omega)\wedge dt + \frac{1}{2}\omega^2$.
\epr
\begin{remark} If $L$ is not a sphere, then the same proof shows that $\bm_{Y,\phi_\lambda}= \bm_L\times f^{-1}(0)$ 
for $\lambda\not=0$ small enough.
This remains an obstructed situation, in the $G_2$ point of view.
\end{remark}

\textbf{With boundary. }
Recall that if $\Sigma$ is a complex surface of $N $ and $t\in S^1$, then $X = \Sigma \times \{t\}$ is a coassociative submanifold of $M$. Consider  the problem of associative deformations of  $Y = L\times \{t\}$  with boundary in $X$:
\begin{thm}\label{SL-boundary} Let $t\in S^1$ and $L$ be a special Lagrangian submanifold in a 6-dimensional Calabi-Yau $N$,
such that $L$ has  boundary in a  complex surface $\Sigma$. Let $Y = L\times \{t\}$ in  $N\times S^1$
and $X = \Sigma \times \{t\}$. 
\begin{enumerate}
\item \label{first}
The moduli space $\bm_{Y, X}$ of associative deformations of $L\times \{t\} $ with boundary in the coassociative $\Sigma\times\{t\}$ can be identified with
the moduli space 
 of special Lagrangian deformations of $L$ with boundary in the fixed $\Sigma$. 
\item \label{second}
If the Ricci curvature of $L$ is positive and if the boundary of $L$ has positive mean curvature in $L$, 
 then $\bm_{Y, X}$ is locally smooth and has dimension $g$, where $g$ is the genus of $\partial L$. 
\end{enumerate}
 \end{thm} 
  Although the moduli space is 
 smooth, its dimension exceeds by one the
 index of the deformation problem, see the beginning
 of the proof of the second assertion. As a consequence, Theorem \ref{coasso}
 shows that generic perturbations of the boundary condition will decrement by one the dimension of the initial moduli space.\\
 
   Note moreover that the deformation theory in \cite{Bu} concerns 
  \textit{minimal Lagrangian} submanifolds with boundary
  in $\Sigma$, a wider class than that of \textit{special Lagrangian} submanifolds of fixed phase.
 \bpr[ of Theorem \ref{SL-boundary} (\ref{first})]
Firstly, if $M$ is equipped with a closed $G_2$-structure $\phi$, note that an associative submanifold $Y$ with boundary
in a coassociative $X$ minimizes the volume in the relative
homology class $[Y]\in H_3(M,X,\Z)$. Indeed, let
$Z$ be any 3-cycle with boundary in $X$, such that $[Z] = [Y]$. 
There is a 4-chain $S$ with boundary in $X$ and $T$ 
a 3-chain in $X$, such that $ Z-Y =  \partial S + T$. 
Since $\phi$ is a calibration, 
$$ Volume (Z)  \geq  \int_Z \phi 
= \int_Y \phi + \int_{\partial S} \phi +  \int_T \phi = \int_Y\phi = Volume (Y)$$
by Stokes and the fact that $\phi$ vanishes on any coassociative submanifold.
By the same arguments as in the closed case, this proves
the identity of the two moduli spaces. 
\epr

\bpr[ of Theorem \ref{SL-boundary} (\ref{second})]
Consider a special Lagrangian $L$ with boundary $\partial L$ 
in a complex surface $\Sigma$. If  $Y = L\times \{t\}$ and $X = \Sigma \times \{t\}$, it is clear that the orthogonal complement $\nu_X$ of $T\partial Y$ 
in $TX$ is equal as a real bundle to $JT\partial L\oplus \{0\}$, and 
$\mu_X$ is the trivial $n\times$-complex line bundle generated by $\partial_t$, where
$n$ is the inward unit normal vector field  of $\partial Y$ in $Y$. 
We begin by computing the index of the boundary problem. 
This is very easy, since $\mu_X$ is trivial, and by Theorem \ref{boundary}, 
we have $\nu_X \cong T\partial L^*$ as $n\times$-bundles. Hence the index equals 
$-c_1(T\partial L) + 1-g = - (2-2g)+1- g = g-1,$
where $g$ is the genus of $\partial L$.
 Now let $\psi =s+ \tau  \frac{\partial }{\partial t}  $
belonging to $\coker (D,\nu_X) = \ker (D, \mu_X)$, where
 $s$ a section of $NL$ and $\tau \in \Gamma(L,\R)$. 
Let $\alpha = -Js^\vee$. By Proposition \ref{harmonic}, 
$\alpha$ is a harmonic 1-form, and $\tau$ is harmonic (note that 
$Y$ is not closed, so $\tau$ may be not constant). By classical results for harmonic 1-forms, we have: 
$$\frac{1}{2}\Delta |\psi|^2 =  \frac{1}{2}\Delta (|\alpha|^2 + |\tau|^2)= |\nabla_L \alpha|^2 + |d\tau|^2+ \frac{1}{2}\text{Ric } (\alpha, \alpha). $$
Integrating on $L\times \{t\}$, we obtain the equivalent of formula (\ref{integration}):
$$ - \int_{\partial Y} \langle \bd_{\mu_X} \psi,\psi\rangle   d\sigma=  \int_Y  |\nabla_L \alpha|^2 + |d\tau|^2+ \frac{1}{2}\text{Ric } (\alpha, \alpha) dy.$$
Lastly, let us compute the eigenvalues of $\bd_{\mu_X}. $
The constant vector $\frac{\partial }{\partial t}$ over $\partial Y$
lies clearly in the kernel of $\bd_{\mu_X}$. 
By Proposition  \ref{bd}, the other eigenvalue of $\bd_{\mu_X}$ 
is $2H$, with eigenspace generated by $n\times \frac{\partial }{\partial t}$. Over $\partial Y$,  $s$ lies in $JTL\cap \mu_X$, hence is proportional
to $n\times \frac{\partial }{\partial t}$. Consequently, 
$\bd_{\mu_X} \psi  =  2 H s$ and  
$$ - \int_{\partial Y} 2H|s|^2 d\sigma=  \int_Y  |\nabla_L \alpha|^2 + |d\tau|^2+ \frac{1}{2}\text{Ric } (\alpha, \alpha) dy.$$
This equation, the  positivity of the Ricci curvature and the positivity of  $H$ show  that $\alpha$ vanishes and $\tau$ is constant. 
So we see that $\dim \coker (D,\nu_X) = 1$, and by the constant rank theorem, $\bm_{Y,X}$ is locally smooth and of dimension
$\dim \ker (D,\nu_X) = g$.  
 \epr
Theorem \ref{SL-boundary} shows an equivalent  result for deformations of special Lagrangian
 submanifold with metric conditions and boundary in a complex surface. Certainly, a direct proof  would be shorter. But it seems to us that our proof has didactic virtues in our context of associative deformations.\\

\noindent
\textbf{A family of examples where $b^3(M) < \dim \coker D$.}
Let $N$ be a projective Calabi-Yau threefold equipped with an ample holomorphic line bundle
$L$, and $N_d $ be the dimension of $\mathbb P H^0(N,L^d)$.  Take $d$ big enough, so that $N_d(N_d-1)/2>b^3(N\times S^1)$
 and choose $C$
a generic complex curve  defined by the intersection of the vanishing locus of two sections of $L^d$.
Then, its moduli space of complex deformations is of dimension  $ N_d(N_d-1)/2$,
so that the dimension of the kernel of the Dirac operator associated to the 
associative $C\times S^1$ is bigger than $b^3(N\times S^1)$.

\section{Appendix }
We will need the following trivial lemma:
\begin{lemma}\label{lemme}
Let $\nabla$ be the Levi-Civita connection on $M$ and $R$ its curvature tensor. 
For any vector fields $w$, $z$, $u$ and $v$ on $M$,
we have 
\beq 
\nabla (u\times v) &=& \nabla u \times  v + u \times \nabla v \\
R(w,z) (u\times v) &=& R(w,z)u\times v + u\times R(w,z)v.
\eeq
If $Y$ is an associative submanifold of $M$ with normal bundle $\nu$,  $u\in \Gamma (Y, TY)$, $v\in \Gamma (Y, TY)$ and $\eta \in \Gamma(Y, \nu)$,
then 
\beq
 \nabla^\top (u\times v) &=& \nabla^\top u \times  v + u \times \nabla^\top v \\
 \nabla^\perp (u\times \eta) &=& \nabla^\top u \times  v + u \times \nabla^\perp v,
\eeq
where $\nabla^\top= \nabla - \nap$ is the orthogonal projection of $\nabla$ to $TY$. 
\end{lemma}
\bpr
Let $x_1, \cdots, x_7$ be normal coordinates on $M$ 
near $x$, and $e_i = \frac{\partial}{\partial x_{i}}$
their derivatives,  orthonormal at $x$.
We have $$u\times v = \sum_i \langle u\times v, e_i\rangle   e_i = \sum_i \phi(u,v,e_i)e_i,$$
so that at $x$, where $\nabla_{e_j} e_i = 0$,
\begin{eqnarray*}
\nabla (u\times v) &= &\sum_i (\nabla \phi (u,v,e_i) + 
\phi (\nabla u,v,e_i) +  \phi ( u,\nabla v,e_i) 
+ \phi ( u,v, \nabla e_i) )e_i\\
&=& \sum_i ( \phi (\nabla u,v,e_i)+ \phi ( u,\nabla v,e_i)) e_i =  
\nabla u \times  v + u \times \nabla v,
\end{eqnarray*}
because $\nabla \phi = 0$. Now if $u$ and $v$ are in $TY$, 
then we get the result after noting that 
$(\nabla u \times v )^\top =   \nabla^\top u \times v$,
because $TY$ is invariant under $\times$. The last relation is implied by $TY\times \nu \subset \nu$
and $\nu\times \nu \subset TY$. The curvature relation 
is easily derived from the definition $R(w,z) = \nabla_w\nabla_z - \nabla_z\nabla_w - \nabla_{[w,z]}$
and the differentiation of the vector product. 
\epr

\subsection{Proof of Lemma \ref{adjoint}}
In this paragraph, we will assume that the ambient manifold $M$ has a torsion-free $G_2$-structure $(\phi,g)$. 
Consider $Y$ an associative submanifold and $\nu$ its normal bundle in $(M,g)$. 
We begin with the classical lemma 
\begin{lemma} 
For a torsion-free structure, the operator $D$ defined in \ref{Operator} is formally self-adjoint, i.e for $s$ and $s' \in \Gamma(Y,\nu)$,
\begin{equation}\label{self}
\int_Y \langle  Ds,s'\rangle   - \langle s,Ds'\rangle   dy = - \int_{\partial Y} \langle  n\times s,s'\rangle   d\sigma, 
\end{equation}
where $d\sigma$ is the volume  induced by the
restriction of $g$ on the boundary, and $n$ is the inward   unit normal vector of $\partial Y$. 
\end{lemma}
\bpr The proof of this lemma is \textit{mutatis mutandis}
the one for the classical Dirac operator, see
Proposition 3.4 in \cite{BoWo} for example. For the reader's convenience
we give a proof of this. 
\beq 
 \langle Ds,s'\rangle   &= & \langle \sum_i e_i\times \nap_i s, s'\rangle   = -  \sum_i\langle \nap_i s,  e_i\times s'\rangle   \\
& = & -  \sum_i d_{e_i} \langle  s,  e_i\times s'\rangle   + \langle s,\nap_i (e_i\times s')\rangle  \\
& = & -  \sum_i d_{e_i} \langle  s,  e_i\times s'\rangle   + \langle s,\nabla_i^\top e_i\times s'+ e_i\times  \nap_i s'\rangle  .
\eeq
By a classical trick, define the vector field $X\in \Gamma (Y,TY)$ 
by $ \langle  X,w\rangle   = - \langle s, w\times s'\rangle   \ \forall w\in TY.$
Note that the  product on the LHS is on $TY$, 
and the one on the RHS is on $\nu$. Now 
$$
 - \sum_i d_{e_i} \langle  s,  e_i\times s'\rangle    = \sum_i d_{e_i} \langle  X,  e_i\rangle   
 =  \sum_i  \langle \nabla^\top_i X,  e_i\rangle   + \langle  X,\nabla_i^\top e_i\rangle  
 =  \sum_i   \text{div } X  -  \langle s,\nabla_i^\top e_i\times s'\rangle  .
$$
By Stokes we get 
$$\int_Y\langle Ds,s'\rangle   dy =  \int_{\partial Y} \langle X,-n\rangle   d\sigma + \int_Y \langle s,Ds'\rangle   dy 
 =  \int_{\partial Y} \langle s,n\times s'\rangle   d\sigma + \int_Y \langle s,Ds'\rangle   dy,
$$
which is what we wanted.\epr
Now, consider $L$ a subbundle of $\nu_{|\partial Y}$ 
of real rank equal to two and invariant under the action of
$n\times$. Let  $s'\in \Gamma (Y,\nu)$ lying in $\coker (D,L)$. 
This means that for every $s\in \Gamma (Y,\nu)$ with $s_{|\partial Y} \in L$,
we have $\displaystyle \int_Y \langle Ds,s'\rangle   dy= 0$. By the former result, 
we see that this is equivalent to 
$$ \int_Y \langle s,Ds'\rangle   + \int_{\partial Y} \langle n\times s,s'\rangle   = 0.$$
This clearly implies that $Ds' = 0$, and $ s'_{|\partial Y} \perp L,$
because $L$ is invariant under the action of $n\times$. 
So  $s' \in \ker (D,L^\perp)$. 
The reverse inclusion holds too by similar reasons.
\subsection{Proof of Proposition \ref{bd}}
\bpr Let $Y$ be an smooth compact associative with boundary, and   $L$ be  a subbundle of $\nu_{|\partial Y}$ invariant under the action of $n\times$.    It is straightforward to check that $\bd_L$ defined in Definition \ref{delta} does not depend
on the chosen orthonormal frame $\{v,w = n\times v\}$. For every $\psi\in \Gamma(\partial Y, L)$ and $f$ a function,
\beq
\bd_L(f\psi) &=& \pi_L (v\times \nabla_w (f\psi) - w\times \nabla_v (f\psi))\\
&= & f\bd_L \psi+ (d_w f)\pi_L ( v\times \psi) - (d_v f)\pi_L ( w\times \psi)  = f\bd_L \psi
\eeq
because $w\times L$ and $v\times L$ are orthogonal to $L$. 
Now, decompose the connexion $\nabla^\top$ on $TY$ as 
$\nabla^\top = \nabla^{\top \partial} + \nabla^{\perp\partial}$
into its two projections along  $T\partial Y$ and along the normal (in $TY$)  $n$-direction. 
For the computations, choose $v$ and $w = n\times v$ 
 the two orthogonal characteristic directions on $T\partial Y$, 
i.e  $\nabla_v^{\top\partial} n = - k_v v$
 and $\nabla_w^{\top\partial}  n =- k_w w$, where $k_v$ and $k_w$ are the two principal curvatures. 
 We have $\nabla^{\perp \partial}_v v  =  k_v n$ and
 $\langle \nabla^{\perp\partial} _w v,n\rangle   = 0$, and the same, \textit{mutatis mutandis}, for $w$. 
Then, for $\psi$ and  $\phi\in  \Gamma(\partial Y, L)$,  using the fact that $T\partial Y\times L $
is orthogonal to $ L$ and Lemma \ref{lemme} below,
\beq 
\langle \bd_L \psi, \phi\rangle   &=& \langle  \nabla^\perp_w (v\times \psi) - (\nabla_w^{\perp\partial} v)\times \psi- 
\nabla_v^\perp (w\times \psi) + (\nabla^{\perp\partial}_v w)\times \psi,\phi\rangle   \\
&=& \langle  \nabla^\perp_w (v\times \psi) - \nabla^\perp_v (w\times \psi) ,\phi\rangle   
= - \langle v\times \psi, \nabla^\perp_w \phi\rangle   +  \langle w\times \psi ,\nabla^\perp_v\phi\rangle   \\
&= & \langle \psi, v\times \nabla^\perp_w \phi - w\times \nabla^\perp_v \phi\rangle   =  \langle \psi, \bd_L \phi\rangle  .
\eeq
To prove that the trace of $\bd_L$ is $2H$, let $e\in L$ be a local unit section of $L$. 
We have $n\times e \in L$ too, and using again Lemma \ref{lemme} and relation  (\ref{permut}),
\beq
\langle \bd_L(n\times e),n\times e\rangle   &=& \langle  v\times ((\nabla^{\top\partial}_w n)\times e) + v \times (n\times \nap_w e), n\times e\rangle   \\
&& - \langle w\times ((\nabla^{\top \partial}_v n)\times e) - w \times (n\times \nabla^\perp_v e), n\times e\rangle   \\
&=& \langle v\times (- k_w w \times e) -  w\times (-k_v v\times  e) , n\times e\rangle  \\
&&+ \langle v \times (n\times \nabla^\perp_w e) - w \times (n\times \nabla^\perp_v e), n\times e\rangle\\
&= & \langle k_w n \times e + k_v n\times e , n\times e\rangle  +
 \langle w\times \nabla^\perp_w e  + v \times  \nabla^\perp_v e, n\times e\rangle.
\eeq
Using again relations (\ref{permut}) and (\ref{permut1}) and the fact that $n\times $ is an isometry on the orthogonal complement
of $n$, we get
\beq 
\langle \bd_L(n\times e),n\times e\rangle  
&=&  k_w + k_v - \langle n\times (w\times \nabla^\perp_w e +v \times  \nabla^\perp_v e), e\rangle  \\
& = & 2H - \langle   v\times  \nabla^\perp_w e - w \times  \nabla^\perp_v e , e\rangle \\
&=& 2H - \langle \bd_L e, e\rangle  .
\eeq
This shows that $\text{trace } \bd_L = 2H$. 
\epr

\subsection{Computation of $D^2$}
\bpr[ of Theorem \ref{D^2}]
We compute $D^2$ at a point $x\in Y$. For this, we choose normal coordinates on $Y$ and $e_i\in \Gamma(Y,TY)$ their associated 
derivatives, orthonormal at $x$. To be explicit, $\nabla^\top e_i = 0$ at $x$.
Let $\psi \in \Gamma (Y,\nu). $     
      \begin{eqnarray*}
          D^2\psi &= & \sum_{i,j} e_i\times \nabla^\perp_i (e_j \times  \nabla^\perp_j \psi) \\  
         & = & \sum_{i,j} e_i\times (e_j \times \nabla^\perp_i \nabla^\perp_j \psi) 
          +\sum_{i,j} e_i\times (\nabla_i^\top e_j \times  \nabla^\perp_j\psi).
          \end{eqnarray*}
        The second sum of the right hand side vanishes, so that using relations (\ref{permut}) and (\ref{permut1}) for the first sum we get 
           \begin{eqnarray*}
       D^2\psi    &=&
          -\sum_i\nabla^\perp_i \nabla^\perp_i \psi - \sum_{i\not= j} 
(e_i \times  e_j) \times \nabla^\perp_{i}\nabla^\perp_{j} \psi\\
          &=&
   \nabla^{\perp *}\nap \psi - 
  \sum_{i\langle  j} (e_i \times  e_j) \times 
(\nabla^\perp_{i}\nabla^\perp_{j} -  \nabla^\perp_{j}\nabla^\perp_{i}) \psi\\
  &=&   \nabla^{\perp *}\nap \psi - 
  \sum_{i\langle  j} (e_i \times  e_j) \times R^\perp(e_i,e_j) \psi.
  \eeq
  Since $(e_i \times  e_j) \times R^\perp(e_i,e_j)$ is symmetric in $i, j$, this is equal to 
 $$  \nabla^{\perp *}\nap \psi - 
  \frac{1}{2}\sum_{i, j} (e_i \times  e_j) \times R^\perp(e_i,e_j) \psi.
$$
The main tool for what follows is the Ricci equation. Let $u$, $v$ be sections of $\Gamma(Y,TY)$
and $\phi$, $\psi$ be elements of  $ \Gamma(Y,\nu)$.  
$$ \langle R^\perp(u,v)\psi,\phi\rangle   =  \langle R(u,v)\psi,\phi\rangle   + \langle  (A_\psi A_\phi- A_\phi A_\psi) u,v\rangle  ,$$
where $A_\phi (u) = A(\phi)(u)  = -\nabla^\top_u \phi.$
Choosing $\eta_1, \cdots, \eta_4$ an
orthonormal basis of $\nu$ at the point $x$, we get 
\beq
-  \frac{1}{2}\sum_{i, j} (e_i \times  e_j) \times R^\perp(e_i,e_j) \psi &=& - \frac{1}{2}\sum_{i, j, k} \langle (e_i \times  e_j) \times R^\perp(e_i,e_j) \psi,\eta_k\rangle  \eta_k\\
&=& \frac{1}{2} \sum_{i, j, k} \langle R^\perp(e_i,e_j) \psi,(e_i \times  e_j) \times \eta_k\rangle  \eta_k\\
&=& -\frac{1}{2} \pi_\nu \sum_{i, j} (e_i \times  e_j) \times R(e_i,e_j) \psi \\
&&+
\frac{1}{2}\sum_{i, j, k}\langle (A_\psi A_{(e_i \times  e_j)\times\eta_k}- A_{(e_i \times  e_j)\times\eta_k} A_\psi) e_i , e_j\rangle  \eta_k.
\eeq
Using the classical Bianchi relation $R(e_i,e_j)\psi = -R(\psi, e_i)e_j - R(e_j,\psi) e_i$,
the first part of the sum $-\frac{1}{2} \pi_\nu \sum_{i, j} (e_i \times  e_j) \times R(e_i,e_j) \psi$ is equal to 
\beq
I = - 2\pi_\nu (e_1\times R(e_2,\psi)e_3 + e_2\times R(e_3,\psi) e_1 + e_3 \times R(e_1,\psi)e_2) =\\
- 2\pi_\nu (e_1\times R(e_2,\psi)(e_1\times e_2) + e_2\times R(e_3,\psi) (e_2\times e_3) + 
e_3 \times R(e_1,\psi)(e_3\times e_1) )=\\
- 2\pi_\nu (e_1\times (R(e_2,\psi)e_1\times e_2+ e_1\times R(e_2,\psi)e_2) + 
                e_2\times (R(e_3,\psi) e_2\times e_3+  e_2\times R(e_3,\psi)e_1)+\\
                e_3\times (R(e_1,\psi) e_3\times e_1+  e_3\times R(e_1,\psi)e_2)) =\\
                - I + 2\pi_\nu \sum_i R(e_i,\psi) e_i,
\eeq
which gives $I = \pi_\nu \sum_i R(e_i,\psi) e_i$.
The Weingarten endomorphisms are symmetric, so that the second part of the  sum is 
\beq 
\frac{1}{2}\sum_{i, j,k}\langle  A_{(e_i \times  e_j)\times\eta_k}e_i,A_\psi e_j\rangle   \eta_k
- 
\frac{1}{2}\sum_{i, j,k}\langle  A_\psi e_i,A_{(e_i \times  e_j)\times\eta_k}e_j\rangle   \eta_k.
\eeq
It is easy to see that the second sum is the opposite of the first one. 
We compute 
\beq 
A_{(e_i \times  e_j)\times\eta_k}e_i &=&-(\nabla^\perp_i e_i \times e_j)\times \eta_k
 - (e_i \times \nabla^\perp_i e_j)\times \eta_k
  + (e_i\times e_j)\times A_{\eta_k} e_i.
\eeq
But we know that an associative submanifold is minimal, so that
$ \sum_i \nabla^\perp_i e_i =0.$
Moreover, differentiating the relation $e_3 = \pm e_1\times e_2$, one 
easily checks that $ \sum_i e_i\times \nabla^\perp_j e_i=0.$
Summing, the only remaining term is 
$$
\sum_{i, j, k}\langle  (e_i\times e_j)\times A_{\eta_k}e_i,A_\psi e_j\rangle   \eta_k. $$
We now use the classical formula for vectors $u$, $v$ and $w$ in $TY$:  
$$ (v\times w) \times u = \langle u,v\rangle  w - \langle u,w\rangle   v, $$
hence $$(e_i\times e_j)\times A_{\eta_k}e_i = \langle A_{\eta_k}e_i , e_i\rangle   e_j - \langle A_{\eta_k}e_i , e_j\rangle   e_i.$$
One more simplification comes from  $\sum_i \langle  A_{\eta_k}e_i,e_i\rangle   =0$ 
for all $k$ because  $Y$ is minimal, 
so our sum is now equal to 
\beq
-\sum_{i, j, k}\langle A_{\eta_k}e_i , e_j\rangle   \langle e_i,  A_\psi e_j\rangle   \eta_k = -\mathcal{A}\psi.
\eeq
\epr

\subsection{Computation of $D$ in the Calabi-Yau extension}
\begin{prf}[ of Proposition \ref{harmonic}] We will use the simple  formula $ \nabla^\perp Js = J\nabla^\top s$ for all sections $s\in \Gamma(L, NL)$.
For $(s,\tau) \in \Gamma(L, NL)\times \Gamma(L,\R)$, and $e_i$ local orthonormal frame on $L$, 
\beq 
D(s,\tau) &=& \sum_{i,j} \langle e_i\times \nabla^\perp_i s,Je_j\rangle  Je_j +   
 \sum_i\langle e_i\times \nabla^\perp_i s, \partial_t\rangle  \partial_t + 
 \sum_i \partial_i \tau \ e_i \times \partial_t\\
     &=& J\sum_{i,j} \phi(e_i,\nabla^\perp_i s, Je_j) e_j 
           + \sum_i \phi(e_i, \nabla^\perp_i s, \partial_t)\partial_t+
           J\sum_{i,j} \partial_i \tau \ \langle e_i\times \partial_t,Je_j\rangle  e_j,
           \eeq
           where we used that $e_i\times \partial_t \perp \partial_t$. 
           \beq
     &=& J\sum_{i,j}\text{Re }\Omega(e_i, \nabla^\perp_i s,Je_j) e_j + 
                  \sum_i \omega(e_i,\nabla^\perp_i s) \partial_t+ 
                J  \sum_{i,j} \partial_i \tau \ \phi(e_i,\partial_t,Je_j)e_j \\
     &= & J\sum_{i,j} \text{Re }\Omega (e_i, J\nabla_i^\top\sigma, Je_j)e_j +
                                \sum_i \omega(e_i, J\nabla_i^\top \sigma) \partial_t +
                          J      \sum_{i,j} \partial_i \tau \ \omega (Je_j, e_i) e_j,
                                \eeq
 where $\sigma = -Js \in \Gamma(L,TL)$.
                                \beq
     &= & -J \sum_{i,j} \text{Re }\Omega (e_i, \nabla_i^\top\sigma, e_j)e_j + 
                                \sum_i \langle e_i,\nabla_i^\top \sigma\rangle   \partial_t -
                        J   \sum_{i,j} \partial_i \tau \langle e_j, e_i\rangle   e_j\\                          
     &= & -J \sum_{i,j} Vol (e_i, \nabla_i^\top\sigma, e_j) e_j + 
                                \sum_i \langle e_i,\nabla_i^\top \sigma\rangle   \partial_t - 
                                                          J \sum_{i} \partial_i \tau e_i  ,                            
                                                             \eeq
    since $\text{Re }\Omega $ is the volume form on $TL$.
It is easy to find that this is equivalent to  
$$D(s,\tau) =  -J ({*}d\sigma^\vee)^\vee + ({*}d*\sigma^\vee) \partial_t - J (d\tau)^\vee,$$
and so $D^\vee (\sigma^\vee, \tau) = (-{*}d \sigma^\vee - d\tau, {*}d*\sigma^\vee).$
Now, since $d^* = (-1)^{3p+1}{*}d*$ on the $p$-forms, one easily checks the formula for $D^2$. 
 \end{prf}


\providecommand{\bysame}{\leavevmode\hbox to3em{\hrulefill}\thinspace}
\providecommand{\MR}{\relax\ifhmode\unskip\space\fi MR }
\providecommand{\MRhref}[2]{%
  \href{http://www.ams.org/mathscinet-getitem?mr=#1}{#2}
}
\providecommand{\href}[2]{#2}

\noindent 
D.~\textsc{Gayet}\\
Universit\'e de Lyon, CNRS, Universit\'e Lyon 1, Institut Camille Jordan,\\
 F--69622 Villeurbanne Cedex, France\\
e-mail: \texttt{gayet@math.univ-lyon1.fr}

\medskip

\end{document}